\documentclass[1p]{elsarticle}

\usepackage[T1]{fontenc}
\usepackage{textcomp}
\usepackage[english]{babel}
\usepackage{hyperref} 

\usepackage{amsmath,amssymb,theorem} 
\usepackage{mathtools}
\usepackage[squaren,Gray]{SIunits}

\usepackage{graphicx}
\usepackage{subcaption}
\usepackage[section]{placeins}

\usepackage{soul}
\usepackage[dvipsnames]{xcolor}

\hypersetup{
  colorlinks   = true, 
  urlcolor     = blue, 
  linkcolor    = blue, 
  citecolor   = red 
}
\renewcommand{\P}{\partial}
\newcommand{\eps}{\ensuremath{\varepsilon} }

\newcommand{\epsOrder}[1]{
    \ifcase#1
    \ensuremath{\mathcal{O}(1) }\or\ensuremath{\mathcal{O}(\eps) }
    \else
    \ensuremath{\mathcal{O}(\eps^{#1}) }
    \fi}

\newcommand{\tagreview}[1]{{ #1} }

\newcommand{\Ct}{\ensuremath{Ct} }
\newcommand{\We}{\ensuremath{\mathrm{We}} }
\newcommand{\Nu}{\ensuremath{\mathrm{Nu}} }
\newcommand{\Pe}{\ensuremath{\mathrm{Pe}} }
\renewcommand{\Re}{\ensuremath{\mathrm{Re}} }
\newcommand{\Bi}{\ensuremath{\mathrm{Bi}} }
\newcommand{\Bit}{\ensuremath{\tilde{\mathrm{Bi}}} }

\newcommand{\Ka}{\ensuremath{\mathrm{Ka}} }

\graphicspath{{figures/}}

\begin{document}

\begin{frontmatter}

	\title{A new family of reduced models for non-isothermal falling films.}
	\author[usmb]{Nicolas Cellier\corref{cor1}}
	\ead{contact@nicolas-cellier.net}
	\author[usmb]{Christian Ruyer-Quil\corref{cor2}}
	\ead{christan.ruyer-quil@univ-smb.fr}

	\begin{keyword}
		heat transfer;
		falling films;
		asymptotic expansion
	\end{keyword}


	\address[usmb]
	{
		LOCIE, UMR 5271, CNRS, Universite Savoie Mont~Blanc, Le Bourget-du-Lac, France
	}

	\begin{abstract}
		New asymptotic models are formulated to capture the thermal transfer across falling films. These models enable us to simulate a wide range of Biot and Peclet number values, without displaying nonphysical behaviors. The models correctly capture the onset of the thermally developed regime at the inlet of the flow. To evaluate the parameter space of acceptability, a comparison has been made with the primitive equation solution for periodic boundary conditions, as well as for an open flow with a periodic forcing at the inlet. A good agreement is obtained for moderate to high Peclet numbers.
	\end{abstract}

\end{frontmatter}

\begin{table}[htbp]
	\centering
		\begin{tabular}{cccccccc}
			\hline
			$x$      & streamwise coord.   & $T_w$             & wall temp.                     & $\We$  & Weber number       \\
			$y$      & cross-stream coord. & $T_a$             & ath. temp.                     & $\Ct$  & Inclination number \\
			$t$      & time                & $h$               & film thickness                 & $\Pe$  & Peclet number      \\
			$g$      & grav. acc.          & $q$               & local flow rate                & $\Bi$  & Film Biot number   \\
			$\beta$  & plate inclination   & $\theta$          & $T|_{y=h}$                     & $\Bit$ & Biot number        \\
			$\nu$    & kinematic visc.     & $\varphi$         & $h^{2}\,\partial_{yy}T|_{y=h}$ & $\Pr$  & Prandtl number     \\
			$\mu$    & dynamic visc.       & $T$               & fluid temp.                    & $\Nu$  & Nusselt number     \\
			$\rho$   & specific mass       & $u$               & velocity (x)                   &        &                    \\
			$\alpha$ & diffusivity         & $v$               & velocity (y)                   &        &                    \\
			$k$      & conductivity        & $X_{\mathrm{Nu}}$ & Nu. flow equiv. var.           &        &                    \\
			$H$      & conv. coeff.        & $\Re$             & Reynolds number                &        &                    \\
			\hline
		\end{tabular}
		\caption{Nomenclature}
		\label{table:nomenclature}
\end{table}

\section*{Introduction}
Falling films form thin layers of liquid flowing on a tilted plate, with a thickness of the order of a millimeter or less. Starting with the works of \citet{nusselt1916}, followed by \citet{kapitza1948}, this topic has been heavily studied and the hydrodynamic of a tilted falling film flowing on a smooth plate is well known. Curiously, the interplay between heat or mass transfers and the wavy dynamics of falling films has been far less studied, even though \citet{frisk1972,yoshimura1996} demonstrated that the wavy regime of the film can indeed increase several folds the heat and mass transfer coefficients between the liquid and the gas. Most studies on heat and mass transfer across the film focus on the wave-less smooth film situation \cite{killion2001}. Only a few studies have been devoted to the wavy regime and addressed heat transfer and hydrodynamics couplings by solving the Fourier equation across the film, the hydrodynamics being dealt with the Navier-Stokes equations or a reduced model. The former approach leads to expensive computations, hardly compatible with a parametric study of the phenomena. It has been restricted either to 2D simulations in a domain corresponding to a full exchanger plate~\cite{serifi2004}, or to numerical domains of limited extensions ~\cite{haroun2012,haroun2010,trifonov2014,nguyen2012,miyara1999}. The latter approach allows better performances with acceptable accuracy but is still not fast enough to allow extensive studies, such as sensibility analyses or optimizations. This explains a lack of numerical studies of non-isothermal falling films, the computation being too expensive to simulate evolution on a full exchanger plate for a significant time interval at a reasonable cost.

Another approach is to use reduced models for both fluid dynamic and heat transfer as proposed by \citet{hirshburg1982} some years ago or more recently by \citet{aktershev2019} but under a fully-linear assumption for the temperature field. The convection effect has thus been neglected. Later, \citet{ruyer-quil2005} have developed a reduced model, based on the weighted residual integrated boundary layer (WRIBL) approach. They considered a linear distribution of the temperature field as a closure hypothesis. This linear distribution was parameterized with only one variable corresponding to the free surface temperature \(\theta = T|_{\bar{y}=1}\). Unfortunately, this model shows nonphysical behaviors at large Peclet numbers, as the temperature field may lie outside the admissible range. \tagreview{\citet{trevelyan2007a} proposed a family of models for both constant temperature and imposed heat flux conditions at the wall. Their approach improves over \citet{ruyer-quil2005} by using test functions verifying all boundary conditions. However, occurrences of non-physical negative temperatures are still observed as the Peclet number is raised.} Recently, \citet{chhay2017} derived a one-equation model in a conservative and Galilean-transform invariant form. In that case, the temperature is bounded, but the model introduced a non-physical critical value \(\theta_c = 7 / 22\) at which convective terms cancel out. It seems that a more complex parameterization of the temperature field is required to overcome this deficiency. \tagreview{Lastly, \citet{thompson2019} considered the interaction of a falling film with a non-uniform heating and derived second-order consistent models. However, the inclusion of second-order convective terms limits their applicability to low and moderate Peclet numbers. This is due to the non-physical vanishing of the diffusion terms at a critical Peclet number. The main focus of our study is to overcome the aforementioned limitations of previous attempts and reach a moderate to high Peclet number domain of applicability with reduced models.}

A cure to these shortcomings have been proposed with a two fields parametrization (\citet{cellier2018b}), adding a second variable \(\phi = \partial_y T\,h|_{\bar{y}=0}\) corresponding to the heat flux at the wall. This approach leads to better results than before but still shows some nonphysical behaviors, with an unbounded temperature. Moreover, the damping terms were not correctly accounted for, non-stationary simulations showing a faster development of the wall thermal boundary layer in the case of the primitive Fourier equation than predicted by the model.
While not being an issue when studying fully developed traveling waves, it impedes the simulation of heat transfer whenever the thermal healing length is not negligible compared to the exchanger dimension. This is especially the case when the Peclet number is increased.

At the same time, only a few experimental studies on this topic are available due to numerous difficulties. Thermocapillarity (Marangoni effect) leads to dry patches formation that is highly problematic when a tracer is used in the fluid. Properly probing the temperature field inside the fluid depth (and not only the surface temperature) is not trivial. Promising approaches involve planar laser-induced fluorescence (PLIF) or laser-induced luminescence to access the temperature volumetric mean \cite{bruchhausen2005,charogiannis2014,xue2018,schagen2007}. Furthermore, most of the study leads to three-dimensional hydrodynamic regimes, for which an extension of the modeling is required in order to perform a proper comparison \cite{chinnov2013,chinnov2017,chinnov2019}. Difficulties also occur with the control of the boundary condition. Imposing a constant temperature or even a constant flux at the wall is experimentally challenging. A similar issue arises with the heat transfer at the free surface, whose precise monitoring requires to account for the development of thermal boundary layers in the gas flow, a difficulty we disregard here using a Newton law of cooling. In order to avoid these difficulties and to validate precisely our modeling attempts, numerical experiments have been used instead. Solutions to the Fourier equation are used as a reference, which allows us to get rid of the experimental difficulties and focus on the modeling of heated falling film at high Peclet number \label{review:exp}.

In this study, we propose a new non-isothermal falling film model in which the temperature field is parameterized with two variables:
\begin{equation}
    \theta = T|_{\bar{y}=1} \qquad \varphi= h^2\frac{\P^2T}{\P y^2}\Big|_{\bar{y}=1}\,.
\end{equation}
This is an attempt to overcome the limitation of the modeling proposed by \citet{ruyer-quil2005} and \citet{chhay2017} by adding more accuracy in the temperature field representation, and obtain a model that possesses coherent damping rates.

The paper is structured as follows: section~\ref{section:primitive-equations} presents the problem to be solved. Our modeling attempt follows in section~\ref{section:development}. These attempts are validated by the Fourier solutions and former attempts by performing two tests. The first one consists of linearizing the equations assuming a non-deformable interface. Construction of large-amplitude nonlinear traveling waves provides the second test.\label{review:twovalid}
We next discuss time-dependent simulations of heat transfer across falling liquid films in extended domains using our model (section~\ref{section:results}). The accuracy and region usefulness of our models are then discussed in the parameter space Biot versus Peclet numbers.

\section{Primitive equations}
\label{section:primitive-equations}

\subsection*{Notations}

We consider a plane making an angle $\beta$ with the horizontal.
We restrict ourselves to the two-dimensional case where the solution
is independent of the span-wise coordinate, and we introduce $x$ and $y$ to refer to the steam-wise and cross-stream coordinates respectively. A film of thickness $h$ flows on a plane maintained at constant temperature $T_w$ and exchanges heat with a cold atmosphere $T_a$ with a constant heat transfer coefficient $H$.

\begin{figure}[ht]
	\centering
	\includegraphics{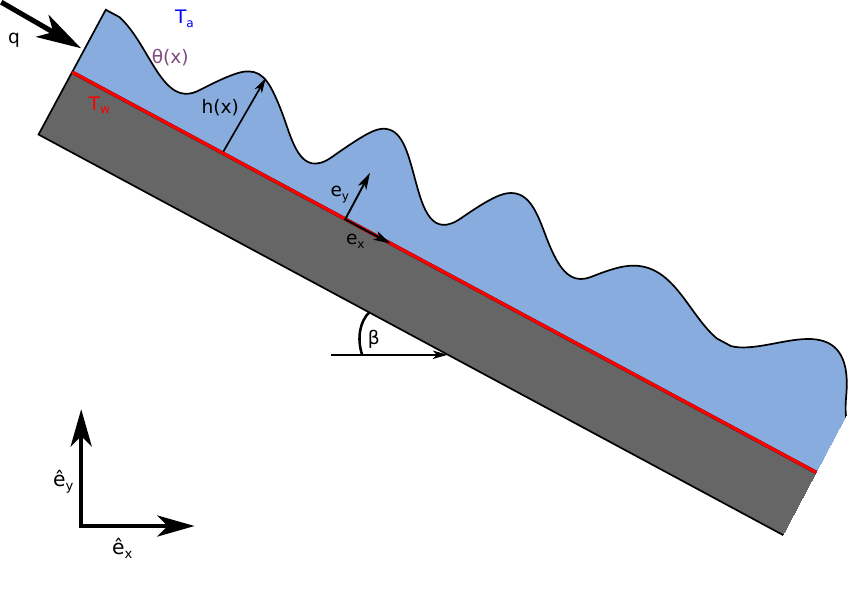}
	\caption{Sketch of a heated falling film (slice).}
\end{figure}

Here we turn directly to dimensionless equations and choose a scaling based on the Nusselt film thickness $h_N = [3 \nu q_L/(g\sin\beta)]^{1/3}$ and the velocity $3 u_N = g\sin \beta h_N^2/\nu$ corresponding to three times the averaged velocity of the Nusselt solution, where $q_L$ is the volumetric flow rate per unit span-wise length, $\nu=\mu/\rho$ is the kinematic viscosity and $g$ is the gravitational acceleration.  Our choice of a velocity scale corresponds to the speed of kinematic waves generated by the deformation of the free surface in the long-wave limit, as traveling waves have a speed close to it.

The dimensionless primitive equations thus consist in the Navier-Stokes equations
\begin{subequations}
\label{eq:primary}
	\begin{eqnarray}
		3 \Re \left(\P_t u + u\P_x u + v \P_y u\right)
		&=& -\P_x p + \P_{xx} u + \P_{yy} u +1\,,\\
		3 \Re \left(\P_t v + u\P_x v + v \P_y v\right)
		&=& -\P_y p + \P_{xx} v + \P_{yy} v\,,\\
		\P_x u + \P_y v &=& 0\,,
	\end{eqnarray}
	the Fourier equation
	\begin{equation}
		\label{eq:fourier}
		3 \Pe\left(\P_t T +  u\P_x T + v \P_y T \right) =
		\P_{xx} T + \P_{yy} T\,,
	\end{equation}
	completed by the no-slip condition at the wall
	\begin{equation}
		u=v=0\qquad \hbox{at}\quad y=0\,,
	\end{equation}
	the kinematic condition at the free surface
	\begin{equation}
		\label{eq:kincond}
		\P_t h + u|_h \P_x h = v|_h\,,
	\end{equation}
	a temperature imposed condition at the wall and  a Newton law of cooling at the free surface
	\begin{align}
		\label{eq:CMA-BC-1}
		T & = 1 \qquad                      & \hbox{at}\qquad y=0\,, \\
		\P_y T - \P_x h \P_x T
		  & = - \Bi T \sqrt{1 + (\P_x h)^2} & \hbox{at}\qquad y=h\,.
		\label{eq:CMA-BC-2}
	\end{align}
	We note that equation \eqref{eq:kincond} is formally equivalent to the mass balance
	\begin{equation}
		\label{eq:mass-balance}
		\P_t h + \P_x q = 0
	\end{equation}
\end{subequations}
where $q=\int_0^h u dy$ is the flow rate.
$\Re = \frac{u_N\,h_N}{\nu} = \frac{q_L}{\nu}$ is the Reynolds number, $\Pe=\Pr\Re$ is the Peclet number and
$\Pr = \frac{\nu}{\alpha}$ is the Prandtl number. Finally, $\Bi = \frac{H\,h_N}{k}$ is the film Biot number, where  $\alpha$, $k$ and $H$ are the thermal diffusivity, the conductivity and the convective heat transfer coefficient. It is also useful to define a second Biot number $\Bit=\frac{H\,l_{\nu}}{k}$ based on a length $l_{\nu}=\left(\frac{\nu^2}{g\,\sin{\beta}}\right)^{1/3}$ corresponding to the balance of gravity and viscosity. In contrast with the film Biot number $\Bi$, the Biot number $\Bit$ is independent of the Reynolds number. The atmosphere has no active effect on the film hydrodynamic and the thermocapillary effect is not taken into account (but can be easily added to the model derivation if needed) as this study focus on developing a model compatible with the high Peclet case. \label{review:hydro-bdc}

\section{Development}
\label{section:development}
In the following, we focus on the derivation of averaged heat equations which enables us to solve the heat transfer within the film more easily than solving the Fourier equation \eqref{eq:fourier} within the framework of the long-wave assumption. We thus introduce a film parameter $\eps$ as the ratio of the typical thickness of the film to the typical length of the waves. The derivatives are of the order of this $\eps$ term, with respect to the stream-wise direction $x$ or with time, as the film evolution is assumed to be slow. As a consequence, the cross-stream velocity is $v = -\int_0^y \partial_x u dy= O(\eps)$.

Within this framework, we further assume that the velocity field $u$ remains close to the parabolic profile corresponding to the Nusselt flow, i.e.
\begin{eqnarray}
	u &=& u^{(0)} + O(\eps)\\
	&=& \frac{3q}{h} \left({\bar{y}} - \frac{1}{2}{\bar{y}}^2\right) + O(\eps)\,,
\end{eqnarray}
where $q=\int_0^h u\,dy$ is the local flow rate and $\bar{y}=y/h$ is a reduced coordinate.
The velocity field is thus parameterized with two variables, the film thickness $h$ and the local flow rate $q$, whose evolution is governed by the mass balance \eqref{eq:mass-balance} and an averaged momentum balance. Several models have been proposed within this framework. Let us cite for instance, the model proposed by Vila and coworkers \cite{boutounet2011}.

\subsection{Gradient expansion approach}
\label{section:GEA}
We thus aim at an integral approximation of the energy balance which mimics the elimination of the cross-stream coordinate that is achieved in Saint-Venant like models. To this aim, we shall project the temperature distribution onto a carefully chosen set of functions. The associated amplitudes of these functions will form our parametrization of the temperature field. The evolution equations associated to these amplitudes will approximate the variations in space and time of the temperature within the waves.
The obtained sets of reduced equations will be validated using two different tests. The first one corresponds to the linear damping eigenmodes of the diffusion operator. As observed in \cite{cellier2017}, passing this test is crucial to correctly capture the thermal entrance region of the film where the thermal boundary layers develop from the wall and free surface.
The second test is the construction of traveling-wave solutions of large amplitude. In the latter case, the thermal regime is developed but differs from the Nusselt linear temperature profile due to convective effects.

Let us first consider that the temperature distribution is never too far from its stationary ($\P_t=0$) and uniform ($\P_x =0$) distribution, i.e. a linear distribution given by:
\begin{equation}
\label{eq:TNu}
	T_{\rm \Nu} = 1 + \left(\frac{1}{1 + \Bi h} -1\right) {\bar{y}}\,.
\end{equation}
A regular expansion around $T_{\rm \Nu}$ with respect of the film parameter $\epsilon$, i.e. $T= T_{\rm Nu}(h) + \epsilon T_1 + \epsilon^2 T_2 \ldots$ is next obtained by solving in sequence the Fourier equation \eqref{eq:fourier} at each order. The result can be found in e.g. \cite{kalliadasis2012}, where the corrections $T_1$, $T_2$ and so on are all functions of $h$ and its derivatives.
Within this framework, the temperature field is thus entirely slaved to the kinematics of the film flow.
However, it is well known that this description of the temperature field is inaccurate
whenever the advection of heat by the flow is non-negligible, i.e. whenever the Peclet number is of order one or larger. We, therefore
revisit the gradient expansion by allowing some degrees of freedom to the temperature distribution. \label{review:Texpand}

Our starting point is the linear relaxation of temperature for a uniform film flow. Considering that the film
thickness $h$ and velocity field ($u$, $v$) are known (and constant),
linearization of the energy balance around the conductive equilibrium and decomposition in normal modes can be done by writing
$T = T_{\rm Nu} + \tilde{T}(\bar{y}) \exp{(i kx + \lambda t)} $, $\tilde{T} \ll T_{\rm Nu}$, where $\lambda$ is the eigenvalue, $k$ a real wavenumber and again ${\bar{y}} = y/h$.
\begin{equation}
\label{eq:heat-transfer-alice}
	3\Pe  h^2(\lambda - i\,k u ) \tilde{T} = \P_{{\bar{y}}{\bar{y}}} \tilde{T} \equiv {\cal L} \tilde{T} \quad
	\hbox{with} \quad \tilde{T}|_{{\bar{y}}=0} = 0
	\quad \hbox{and} \quad \P_{\bar{y}} \tilde{T} |_1 + \Bi h \tilde{T} |_1 =0
\end{equation}
Solutions to \eqref{eq:heat-transfer-alice} form discrete branches, as setting $k$ to zero (very long-wave limit) yields
 eigenfunctions $v_{k=0}^{(n)}({\bar{y}})$ and eigenvalues $\lambda_{k=0}^{(n)}$ given by
\begin{subequations}
	\label{eq:lin}
	\begin{equation}
		v_{k=0}^{(n)} = \sin(l_n\,\bar{y}), \quad \lambda_{k=0}^{(n)} = - \frac{l_n^2}{3\Pe}
	\end{equation}
	where $l_n$ are solutions to
	\begin{equation}
		\label{eq:def-l}
		l \cot l +\Bi h  = 0\,.
	\end{equation}
\end{subequations}
All eigenvalues $\lambda^{(n)}$ have a negative real part. They correspond to
relaxation modes promoted by the diffusion of heat across the film.
Two limits are worth investigating. The first one is $\Bi=0$ which
corresponds to an insulated free surface, in which case the $T_{\rm Nu}= 1$ is constant and equal to its value at the wall. The second limit is
$\Bi\to\infty$, in which case the free surface is at the constant
temperature $T=0$ (equal to the temperature of the gas phase).

For $\Bi=0$,the discrete spectrum of ${\cal L}$ for $k=0$ is
\begin{subequations}
	\label{eq:lin-B0}
	\begin{eqnarray}
		l_1 &=& \frac{\pi}{2}\,,\qquad
		l_2 = \frac{3\pi}{2}\,,\qquad
		l_3 = \frac{5\pi}{2}\,,\\
		3\Pe\lambda_{k=0}^{(1)} &\approx& -2.47\,,\qquad
		3\Pe\lambda_ {k=0}^{(2)}\approx -22.21\,,\qquad
		3\Pe\lambda_{k=0}^{(3)} \approx - 61.69\,.
	\end{eqnarray}
\end{subequations}
For $\Bi\to\infty$ we have instead
\begin{subequations}
	\label{eq:lin-Binf}
	\begin{eqnarray}
		l_1 &=& \pi\,,\qquad
		l_2 = 2\pi\,,\qquad
		l_3 = 3\pi\,,\\
		3\Pe\lambda_{k=0}^{(1)} &\approx& -9.87\,,\qquad
		3\Pe\lambda_{k=0}^{(2)} \approx -39.48\,,\qquad
		3\Pe\lambda_{k=0}^{(3)} \approx - 88.83\,.
	\end{eqnarray}
\end{subequations}
Considering long-time evolutions of the temperature, deviations from the linear temperature distribution \eqref{eq:TNu} are all damped by the relaxation eigenmodes. As a consequence, the temperature field is slaved to the film thickness. At shorter time scales,
only the eigenmodes with sufficiently small eigenvalues are effective and the first eigenmodes \eqref{eq:lin} should be taken into account, in which case the temperature field depends not only on $h$ but also on the amplitudes of some eigenmodes. Roberts \cite{roberts2014}
used the center manifold approach to extend this idea in the case of non-uniform film thickness and large deviations. \label{review:manifold}
Following Roberts, we shall assume that the time evolution of the temperature is determined by the evolution on a manifold that is tangent to the first eigenmodes \eqref{eq:lin}.

Let us thus decompose the temperature field into
\begin{equation}
	\label{eq:decomp}
	T = T_{\rm Nu} + T^{(0)}+ T^{(1)}\,,
\end{equation}
where $T^{(0)}$ is aligned with the two first eigenmodes $v_{k=0}^{(1)}$ and $v_{k=0}^{(2)}$. This idea is similar to the semi-analytical method for solving the problem of heating of a uniform film flow that has been proposed by Aktershev and Bartashevich \cite{aktershev2017}.
However, instead of projecting the temperature field on the sinus functions $v_{k=0}^{(n)}$, as proposed by Aktershev and Bartashevich, where $l_n$ are not given explicitly but indirectly through the solution to \eqref{eq:def-l}, it is more convenient to use polynomial approximations.

\label{review:polynomial}
Requiring that ${\tilde v_1(\bar{y})}$ and ${\tilde v_2(\bar{y})}$ are polynomials of the lowest degrees in $\bar{y}$ and $\Bi$  which verify
\begin{subequations}
	\begin{eqnarray}
		{\tilde v}_1'(1) + \Bi h {\tilde v}_1(1) &=&0 \,,\qquad {\tilde v}_1(0)=0\,,\\
		\hbox{for}\qquad \Bi h=0\,,\qquad
		{\tilde v}_1(1) &=&1\,,\qquad {\tilde v}_1'(1) =0\,,\\
		\hbox{for}\qquad \Bi h \gg 1\,,\qquad
		{\tilde v}_1(1) &\ll&\Bi h \,,\qquad {\tilde v}_1'(1)-1 \ll \Bi h \,,
	\end{eqnarray}
\end{subequations}
and
\begin{subequations}
	\begin{eqnarray}
		{\tilde v}_2'(1) + \Bi h {\tilde v}_2(1) &=&0 \,,\qquad {\tilde v}_2(0)=0\,,\\
		\hbox{for}\qquad \Bi h=0\,,\qquad
		{\tilde v}_2(1) &=&1\,,\qquad {\tilde v}_2'(1) =0\,,\qquad {\tilde v}_2(2/3) =0 \,,\\ \nonumber
		\hbox{for}\qquad \Bi h \gg 1\,,\qquad
		{\tilde v}_2(1) &\ll&  \Bi h  \,,\qquad {\tilde v}_2'(1) -1 \ll \Bi h\,, \\
		{\tilde v}_2(1/2) &\ll& \Bi h\,,
	\end{eqnarray}
\end{subequations}
then gives
\begin{eqnarray}
	{\tilde v_1} &=& {\bar{y}}(2-{\bar{y}}) + \Bi h\, {\bar{y}}(1-{\bar{y}})\,,\\
	{\tilde v_2} &=& -12{\bar{y}}\left(\frac{2}{3}-{\bar{y}}\right)
	\left(\frac{5}{4}-{\bar{y}}\right)+ \Bi h\,  2{\bar{y}} (1 - {\bar{y}})\left({\bar{y}}-\frac{1}{2}\right)\,.
\end{eqnarray}
${\tilde v}_1$ and ${\tilde v}_2$ are polynomial approximations to the relaxation
eigenmodes $v_{k=0}^{(1)}$ and $v_{k=0}^{(2)}$. Obviously, these approximations are more accurate at low values of the $\Bi$ number than at high values. Indeed, we anticipate that the most challenging phenomenon to capture is the onset of thermal boundary layers in the vicinity of the hyperbolic stagnation points appearing with the recirculation zones in large-amplitude solitary waves \cite{trevelyan2007a}. These thermal boundary layers do not develop in the limit of large $\Bi$ numbers as the free surface temperature becomes constant and we therefore focus on accuracy on low or moderate values of $\Bi$.
We next introduce a linear combination of ${\tilde v}_1$,  ${\tilde v}_2$ and two variables to represent the departure of the temperature field from the linear temperature distribution. The choice of these variables is particularly important. In order to fully capture the onset of a thermal boundary layer close to the stagnation point at the front of the waves, we chose variables which monitor the temperature distribution close to the free surface.
 The free-surface temperature $\theta=T(y=h)$ is the most obvious choice. We complete it using the derivative of lowest order which is independent to $\theta$.
 As the Newton law \eqref{eq:CMA-BC-2} relates the gradient of temperature to the free-surface temperature, we chose
 $\varphi = h^2\P_{yy} T(y=h)$ such that $\varphi$ has the dimension of a temperature. We thus introduce ${\hat v}_1$ and  ${\hat v}_2$ :
 \begin{subequations}
     \label{eq:hatv1-2}
 \begin{eqnarray}
	\label{eq:hatv1}
 {\hat v_1} &=& {\bar{y}} [3-3{\bar{y}} +{\bar{y}}^2 + \Bi h (2-3{\bar{y}}+{\bar{y}}^2)]\,,\\
 	\label{eq:hatv2}
 {\hat v_2} &=& \frac{1}{2}{\bar{y}}
 \left(1 -{\bar{y}}\right)^2\,.
 \end{eqnarray}
 \end{subequations}
as linear combinations of ${\tilde v}_1$,  ${\tilde v}_2$ verifying
 \begin{subequations}
     \label{eq:hatBC1-2}
 \begin{eqnarray}
	\label{eq:haBC1}
 {\hat v_1}(1) &=& 1, \qquad {\hat v_1}"(1) = 0\,,\\
 	\label{eq:haBC2}
 {\hat v_2}(1) &=& 0, \qquad {\hat v_1}"(1) = 1\,.
 \end{eqnarray}
 \end{subequations}

We first introduce the ansatz
\begin{equation}
	\label{eq:ansatz-1}
	T^{(0)} = [\theta - \theta_0(h)]{\hat v_1}({\bar{y}})
	\qquad \hbox{with} \qquad
	\theta_0 = \frac{1}{1+ \Bi h}\,,
\end{equation}
so that $\theta= (T_{\rm Nu} + T^{(0)})(y=h)$. We emphasize that $T_{\rm Nu} + T^{(0)})(y=h)$ defined by \eqref{eq:ansatz-1} verifies the boundary conditions \eqref{eq:CMA-BC-1} and \eqref{eq:CMA-BC-2}.
Thus, according to our choice of variables, the decomposition \eqref{eq:decomp} with
\eqref{eq:ansatz-1} is set unique by the gauge condition
\begin{equation}
	\label{eq:gauge-1}
	T^{(1)}|_{y=h} = 0\,.
\end{equation}
Inserting the decomposition \eqref{eq:decomp}, \eqref{eq:ansatz-1} in \eqref{eq:fourier} gives
\begin{subequations}
	\label{eq:T2}
	\begin{equation}
		\label{eq:T2-bulk}
		\P_{yy} T^{(1)} = 3\Pe \left(\P_t + u \P_x + v \P_y \right) (T_{\rm Nu} + T^{(0)})
		- \P_{xx}  (T_{\rm Nu} + T^{(0)})
		- \P_{yy}  T^{(0)}
		\,.
	\end{equation}
	Here, the second-order corrections to the advection terms $3\Pe \left(\P_t + u \P_x + v \P_y \right) T^{(1)}$ have been dropped out	while the diffusion terms have been retained. This is justified considering that (i) these corrections are small compared to the other advection terms, (ii) all leading-order physical contributions have been retained in \eqref{eq:T2-bulk}.\tagreview{This procedure is similar to the treatment of the momentum balance using the weighted residuals technique, where second-order inertial terms are dropped from the averaged momentum balance. See the discussion in \citet{richard2016}. Inclusion of these second-order inertial terms is possible but at the expense of complicated formulations or limited ranges of applicability as performed in \citet{scheid2006} using a Padé approximant technique}.\label{review-2ndorder} We note that solving \eqref{eq:T2-bulk} is similar to looking for $T^{(1)}$ in terms of an expansion $T^{(1)}= \eps T^{(1)}_1 + \eps^2 T^{(1)}_2 + \ldots$ with respect to the film parameter where only the leading-order contributions are retained. Equation \eqref{eq:T2-bulk} is completed with the boundary conditions
	\begin{eqnarray}
		\label{eq:CMA-BC-T2}
		T^{(1)} &=& 0 \qquad \hbox{at}\qquad y=0\,,\\
		\P_y T^{(1)} &=&  \P_x h \P_x (T_{\rm Nu} + T^{(0)})
		- \Bi (T_{\rm Nu} + T^{(0)}) \frac{1}{2} (\P_x h)^2
		\quad \hbox{at}\quad y=h\,.
	\end{eqnarray}
\end{subequations}

Solving \eqref{eq:T2} gives the correction $T^{(1)}$ as a polynomial in $y$ whose coefficients are dependent on the variables $h$, $q$, $\theta$, $\varphi$ and their derivatives. The gauge condition \eqref{eq:gauge-1} then provides an evolution equation for the variables $\theta$, namely
\begin{eqnarray}
	\nonumber
	3 \Pe \P_t \theta &=&
	3 \Pe \left\{-\frac{3(82 + 19 \Bi h)}{7(27+ 7 \Bi h)}\frac{q}{h}\P_x \theta
	- \frac{57 \Bi h}{7(27+ 7 \Bi h)} \frac{q\theta}{h^2}\P_x h
	\right. \\ && \left.
	+  \frac{3 [11+ (-11+  38\Bi h)\theta]}{14(27+ 7 \Bi h) h} \P_x q \right\}
	- \frac{60(1+ \Bi h)}{27 +  7 \Bi h} \frac{\theta - \theta^{(0)}}{h^2}
	\nonumber \\ &&
	+ \P_{xx} \theta + \left( \frac{6 + 3(-2 + 7\Bi h)\theta}{27+ 7 \Bi h} \right) \frac{\P_{xx} h}{h}
	+ \left(\frac{6 + 6(-1 + 2\Bi h)\theta}{27+ 7 \Bi h} \right)\frac{(\P_{x} h)^2}{h^2}
	\nonumber \\ &&
	+ \left( \frac{6(8 + 7\Bi h)\theta}{27+ 7 \Bi h} \right) \frac{\P_{x} h \P_x \theta}{h}\,,
	\label{eq:model-theta}
\end{eqnarray}
referred hereinafter as the $\theta$~model.

In essence, equation \eqref{eq:model-theta} is an averaged energy balance which must be contrasted to the model derived by Ruyer-Quil et al. \cite{ruyer-quil2005} using the method of weighted residuals:
\begin{eqnarray}
	\nonumber
	3 \Pe \P_t \theta &=&
	3 \Pe \left\{-\frac{27}{20}\frac{q}{h}\P_x \theta
	+  \frac{7}{40}\frac{(1-\theta)}{h} \P_x q \right\}
	- 3  \frac{\theta - \theta^{(0)}}{h^2}
	\\ &&
	+ \P_{xx} \theta + (1- \theta) \frac{\P_{xx} h}{h}
	+ \left(1-\theta -\frac{3}{2}\Bi h \right)\frac{(\P_{x} h)^2}{h^2}
	+ \frac{\P_x h \P_x \theta}{h}
	\label{eq:EI-Scheid}
\end{eqnarray}
The two energy balances are consistent with the long-wave expansion up to first-order for the convection terms and second-order for the diffusion terms. In fact, \eqref{eq:EI-Scheid} can be obtained following our approach with the ansatz
\begin{equation}
	\label{eq:ansatz-Scheid}
	T^{(0)} = [\theta - \theta_0(h)] {\bar y} \quad \hbox{so that}\quad T_{\rm Nu} + T^{(0)} = 1+ (\theta-1) {\bar y}
\end{equation}

\begin{figure}
	\begin{subfigure}[t]{0.49\textwidth}
		\includegraphics{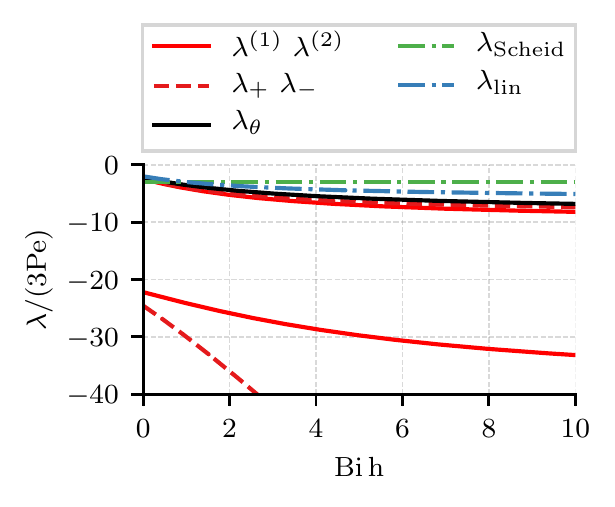}
		\caption{$k=0$}
		\label{fig:lambda-1}
	\end{subfigure}~
	\begin{subfigure}[t]{0.49\textwidth}
		\includegraphics{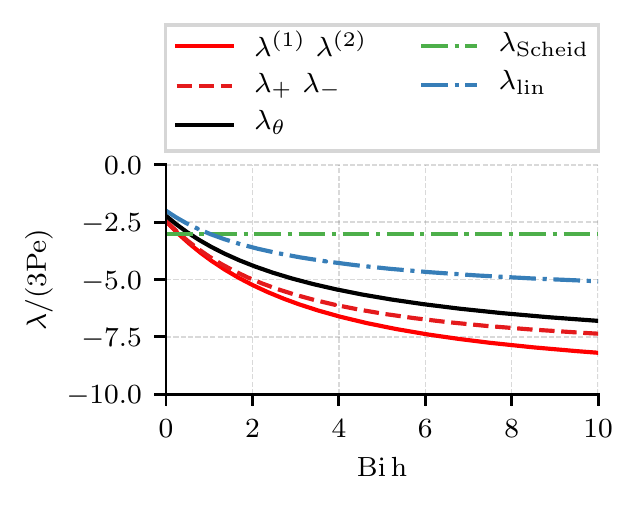}
		\caption{$k=0 $ (enlargement of panel~a)}
		\label{fig:lambda}
	\end{subfigure}
	\begin{subfigure}[t]{0.99\textwidth}
		\includegraphics{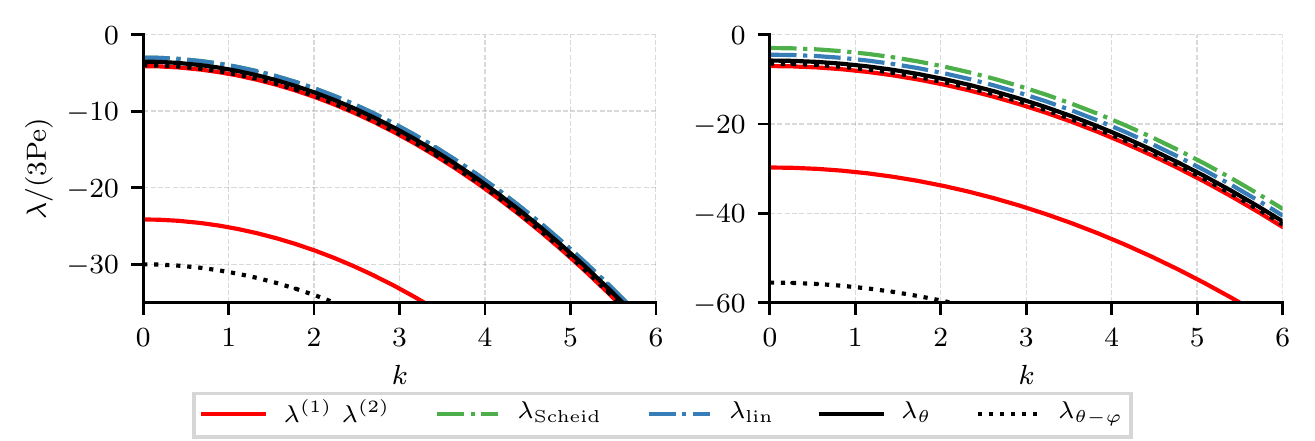}
		\caption{Left, $\Bi=1$ ; right, $\Bi=5$.}
		\label{fig:lambda_k}
	\end{subfigure}
	\caption{Real part of the eigenvalues corresponding to the modal response of the temperature field to a perturbation of wavenumber $k$. Solutions $\lambda_\theta$, $\lambda_\pm$, $\lambda_{\rm Scheid}$ to the models \eqref{eq:model-theta}, \eqref{eq:model-theta-phi}, \eqref{eq:EI-Scheid} and \eqref{eq:model-alice} are compared to the solutions $\lambda^{(n)}$ to \eqref{eq:def-l}}
\end{figure}

Considering the aforementioned first test consisting of the linear relaxation of the temperature to the linear Nusselt distribution, \eqref{eq:model-theta} represents a significant improvement over the former formulation \eqref{eq:EI-Scheid}. Linearizing \eqref{eq:model-theta}
considering a flat film (i.e. $h=1$ and $q=1$) gives a damping rate
\begin{equation}
	\label{eq:l-theta}
	3\Pe\lambda_\theta = - \frac{60(1+ \Bi h)}{27 +  7 \Bi h} - k^2
\end{equation}
which is compared to the eigenvalue
$\lambda^{(1)}$ in figure~\ref{fig:lambda}.
As expected, a much better agreement is observed with the new formulation \eqref{eq:model-theta} than with \eqref{eq:EI-Scheid}.

We note that the projection approach followed by \citet{thompson2019} corresponds to the ansatz \eqref{eq:ansatz-1} where $T^{(0)}$ remains aligned with the Nusselt linear temperature distribution, i.e.
\begin{equation}
	\label{eq:ansatz-Alice}
	T^{(0)} = [\theta - \theta_0(h)] {\hat v}_{\rm lin}
	\qquad \hbox{with} \qquad
	{\hat v}_{\rm lin} = T_{\rm Nu}/\theta_0(h)
\end{equation}
so that ${\hat v}_{\rm lin}(1)=1$ as required by the definition of $\theta$.
	Following our approach, the resulting evolution equation for the
	free-surface temperature then reads
\begin{eqnarray}
	\nonumber
	3 \Pe \P_t \theta &=&
	3 \Pe \left\{-\frac{3(25 + 7\Bi h)}{20(3+ \Bi h)}\frac{q}{h}\P_x \theta
	- \frac{21 \Bi h}{20(3+ \Bi h)} \frac{q\theta}{h^2}\P_x h
	\right. \\ && \left.
	+  \frac{27 \Bi\theta}{20(3+ \Bi h)} \P_x q \right\}
	- \frac{6(1+ \Bi h)}{3 + \Bi h} \frac{\theta - \theta^{(0)}}{h^2}
	\nonumber \\ &&
	+ \P_{xx} \theta + \left( \frac{3\Bi h \theta}{3+ \Bi h} \right) \frac{\P_{xx} h}{h}
	+ \left(\frac{3\Bi h\theta}{3+ \Bi h} \right)\frac{(\P_{x} h)^2}{h^2}
	\nonumber \\ &&
	+ \left( \frac{6(1 +\Bi h)\theta}{3+ \Bi h} \right) \frac{\P_{x} h \P_x \theta}{h}\,,
	\label{eq:model-alice}
\end{eqnarray}
\tagreview{Equation~\eqref{eq:model-alice} represents a truncation of the model derived by \citet{thompson2019} (equation~6.6 in this publication) by dropping second order convective terms proportional to $\Pe^2$.} The corresponding damping rate
\begin{equation}
	3\Pe\lambda_{\rm lin} = - \frac{6(1+ \Bi h)}{3 + \Bi h} - k^2
\end{equation}
is again compared to the eigenvalue
$\lambda_1$ in figure~\ref{fig:lambda-1}.

\label{review:trev}
Our comparisons to previous attempts of one-variable modeling would not be complete without mentioning the work by Trevelyan et al. \cite{trevelyan2007a}. These authors have constructed Galerkin projections of the temperature field which, in contrast to \cite{ruyer-quil2005}, verify the boundary conditions \eqref{eq:CMA-BC-1} and \eqref{eq:CMA-BC-2}. However, by following strictly the Galerkin approach, their one-variable model, referred to as GST[1] in their work, is not consistent with the long-wave expansion (consistency is however recovered
when the number of variables is larger than three).
Considering GST[1], Trevelyan's choice of polynomial projection is equivalent to
\begin{equation}
	\label{eq:ansatz-Trevelyan}
	T^{(0)} = [\theta - \theta_0(h)] {\hat v}_{\rm Tre}(y/h)
	\qquad \hbox{with} \qquad
	{\hat v}_{\rm Tre}({\bar y}) = \frac{1}{2}{\bar y}
	\left(3 -{\bar y}^2 + \Bi h (1-{\bar y}^2) \right)\,.
\end{equation}
This choice of polynomial profile stems from the requirement that $\P_{yy} T=0$ at the wall, as can be proved easily by writing the  Fourier equation \eqref{eq:fourier}. A consistent evolution equation for $\theta$ similar to \eqref{eq:model-theta} and \eqref{eq:model-alice} can easily be formed from the ansatz \eqref{eq:ansatz-Trevelyan} following the approach developed above. For the sake of brevity, we refrain from writing it.

\subsection{Construction of traveling-wave solutions}
A second validation of the modeling approach is offered by the construction of the traveling-wave solutions to the models. Considering a stationary solution in a frame of reference  $\xi = x - c\,t$, moving at constant speed $c$, the set of partial differential equations reduces to ordinary differential equations which is then recast into an autonomous dynamical system \cite{kalliadasis2012}. This dynamical system of finite dimension is solved using a continuation method by Auto07p software \cite{doedel2007}.
We have focused on solitary-wave solutions to \eqref{eq:model-theta} where the hydrodynamics of the film is modeled by the Saint-Venant equations derived by Vila and coworkers
\begin{subequations}
\label{eq:Vila}
\begin{eqnarray}
	\P_t h &=& - \P_x q\,,\\
	3\Re\P_t q &=& - 3\Re\P_x \left(\frac{q^2}{h} +
	\frac{2}{225} h^5\right) = h -3\frac{q}{h^2} + \We \P_{xxx} h\,.
\end{eqnarray}
\end{subequations}
In this section, we compare the solutions to the different one-variable models of heat transfer to the solutions to the primary Fourier problem \eqref{eq:primary} that we have obtained using a classical pseudo-spectral method (see section~\ref{sec:appendix2} for details).

\begin{figure}
	\begin{center}
		\includegraphics{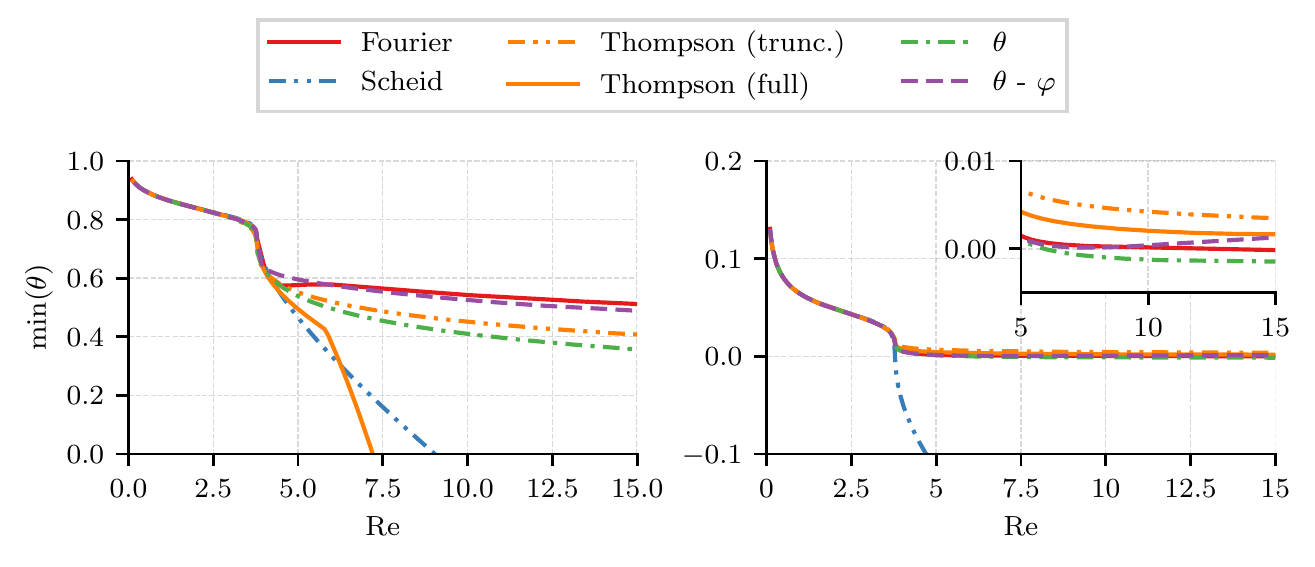}
	\end{center}
	\caption{Minimum of the free surface temperature $\theta$ under a solitary wave as function of the Reynolds number for a vertical water film ($\beta=90^\circ$ and $\Ka=3000$). Left: $\Pr=30$ and $\Bit=0.1$ ; right: $\Pr=7$ and $\Bit=10$. \tagreview{The insert shows an enlargement of the main plot.}
		\label{fig:TW-ts}}
\end{figure}

Equation~\eqref{eq:EI-Scheid} has been shown to be limited to low Peclet values as its solutions present nonphysical values of the free-surface temperature, i.e. $\theta$ lies out of the admissible interval [0,1].
We present therefore in figure~\ref{fig:TW-ts} the evolution of the minimum of
\(\theta\) as a function of the Reynolds number \(\Re\). The film is vertical (\(\beta=90\degree\)) and the liquid properties correspond to water (\(\Ka=3000\)).

Comparisons to the solutions to the Fourier equations \eqref{eq:primary} show that \eqref{eq:model-theta} achieves a much better agreement to the reference than the former averaged energy equation \eqref{eq:EI-Scheid}. Aberrant values of $\theta$ are almost unobserved with the new formulation. Solutions the model by \citet{thompson2019} present a similar property\tagreview{, but only if the second-order convective terms proportional to $\Pe^2$ are dropped out leading to \eqref{eq:model-alice}}. Considering that \eqref{eq:model-theta}, \eqref{eq:EI-Scheid} and \eqref{eq:model-alice} present similar mathematical structures, the origin of the differences in behavior of their solutions is not obvious. In particular, as is generally the case with asymptotic expansions, the conservative structure of the basic Fourier equation is lost. Yet, it can be noticed that among the ansatze for the temperature profile presented so far, \eqref{eq:ansatz-1} and \eqref{eq:ansatz-Alice} verify the Newton law of cooling \eqref{eq:CMA-BC-2}, whereas \eqref{eq:ansatz-Scheid} does not.
In fact, with \eqref{eq:ansatz-Scheid}, the reconstructed temperature field $T=T_{\rm Nu} + T^{(0)} + T^{(1)}$ verifies the Fourier equation, the boundary condition \eqref{eq:haBC1} at the wall, but complies with the Newton law \eqref{eq:CMA-BC-2} at the free surface only in the asymptotic limit where $\theta$ remains close to $\theta_0(h)$, which is a more restrictive condition than ensuring that the second-order advection terms $3\Pe \left(\P_t + u \P_x + v \P_y \right) T^{(1)}$ remains small in comparison to first-order ones.
However, starting with the temperature ansatz \eqref{eq:ansatz-Trevelyan} corresponding to the work by Trevelyan et al. \cite{trevelyan2007a}, the obtained model does present occurrences of negative temperature even though \eqref{eq:ansatz-Trevelyan} verifies the boundary condition \eqref{eq:CMA-BC-2}. Therefore, we conclude that requiring that the temperature profile verifies the
boundary conditions is not sufficient to guarantee that $\theta$ remains within the physical range.
To conclude this comparison of our approach with previous attempts, we have added to
figure~\ref{fig:TW-ts} the curves corresponding to the model by \citet{thompson2019} including second-order convective terms. Besides complicating the problem to solve, inclusion of these second-order terms leads to non-physical values of the temperature as the Peclet number is raised, which severely reduces the parameter range for which this model may be useful.

\begin{figure}[hbtp]
	\begin{center}
		\includegraphics{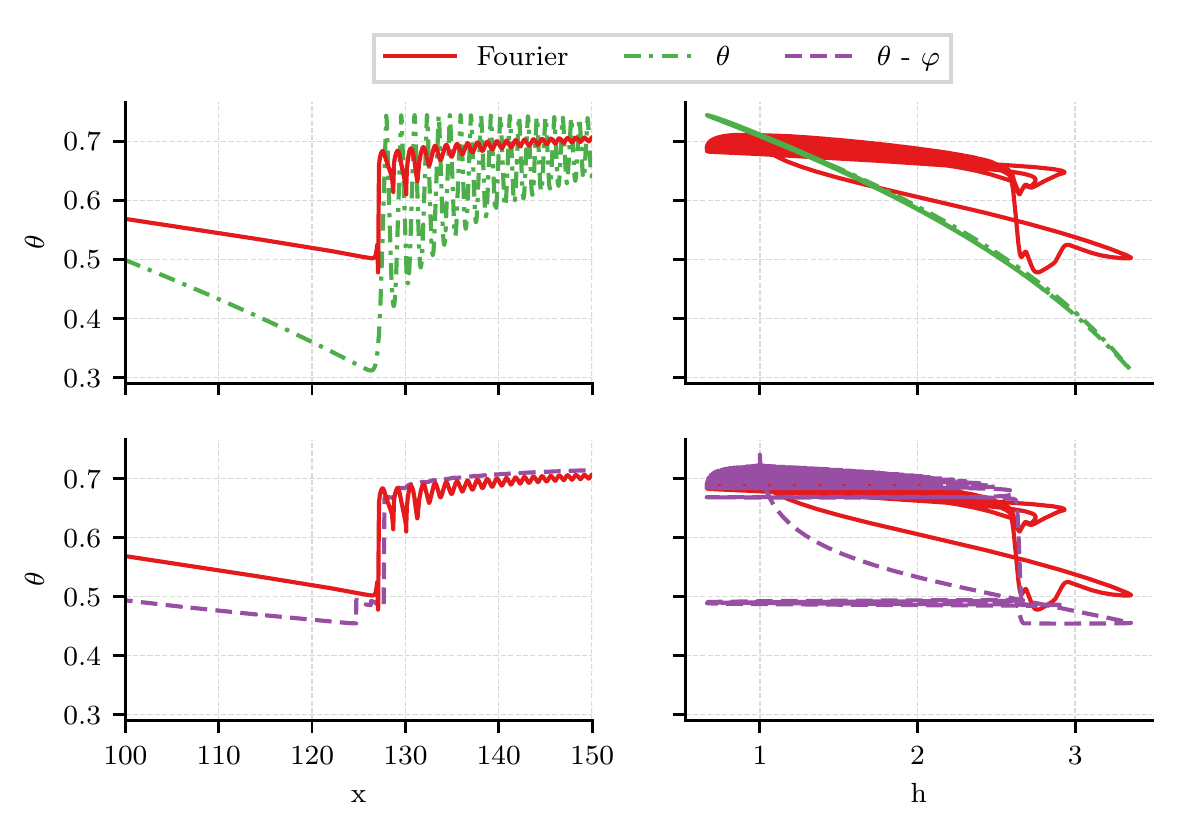}
	\end{center}
	\caption{Distribution of the free surface temperature $\theta$ under a solitary wave as a function of $x$ (left) and the fluid thickness $h$ (right) for a vertical water film ($\beta=90^\circ$, $\Ka=3000$, $\Re=33$, $\Pr=30$ and $\Bit=0.1$).
		\label{fig:theta_thetah}}
\end{figure}

Yet, a close examination of the distribution of $\theta$ (figure~\ref{fig:theta_thetah}) under a wave shows that the model \eqref{eq:model-theta} overestimates the variations of temperature under the wave. The model also fails to reproduce the jump of free-surface temperature at the front of the wave which is promoted by the presence of a roll in the wave crest. This rapid variation signals the development of a thermal boundary layer in the vicinity of a hyperbolic stagnation point at the front of the crest (at $h\approx 2.6$ for the discussed solitary wave). The onset of a thermal boundary layer cannot be captured by \eqref{eq:EI-Scheid} as $\theta$ tends to be a function of $h$ in that case whenever the Peclet number $\Pe$ is large as observed in figure~\ref{fig:TW-ts}.

\subsection{A two-variable model of heat transfer}
Overcoming the limitations of one-variable averaged heat equations demands to enrich the modeling. We can do so
by selecting another variable which reflects the complexity of the temperature field in the vicinity of the free surface, we thus introduce $\varphi = h^2\P_{yy} T(y=h)$ such that $\varphi$ has the dimension of a temperature.

We then introduce a more complete ansatz
\begin{equation}
	\label{eq:ansatz}
	T^{(0)} = [\theta - \theta_0(h)]{\hat v_1}({\bar{y}}) + \varphi {\hat v_2}({\bar{y}})
\end{equation}
The decomposition \eqref{eq:ansatz} is made unique by adding
\begin{equation}
	\label{eq:gauge}
	\P_{yy} T^{(1)}|_{y=h} =0\,.
\end{equation}
to the gauge condition \eqref{eq:gauge-1}.
Solving \eqref{eq:T2} then provides an expression of the correction
$T^{(1)}$ that is consistent with the ansatz \eqref{eq:ansatz} and the long-wave expansion up to first-order for the convective terms and second-order for the diffusion ones. The gauge conditions \eqref{eq:gauge} and  \eqref{eq:gauge-1} then yields coupled evolution equations for the variables $\theta$ and $\varphi$, namely
\begin{subequations}
	\label{eq:model-theta-phi}
	\begin{equation}
		\label{eq:model-ts}
		3 \Pe  (\P_t + u|_{y=h} \P_x)\theta =
		\frac{\varphi}{h^2}  +
		2\Bi \P_x h \P_x \theta + \frac{\varphi}{h^2} (\P_xh)^2 + \Bi \theta\P_{xx} h
		+ \P_{xx} \theta
	\end{equation}
	with $u|_{y=h}=\frac{3}{2}q/h$,
	and
	\begin{eqnarray}
		\nonumber
		3 \Pe \P_t \varphi &=& - 3 \Pe \left(\frac{15}{14} \frac{q}{h}\P_x \varphi + E_\varphi\frac{q}{h} \P_x \theta + F_\varphi  \frac{\P_x q}{h}
		+ G_\varphi  \frac{q \theta}{h^2}\P_x h   \right)
		\\ \nonumber && + \frac{1}{h^2}\left\{-60(1+\Bi h) (\theta -\theta_0(h)) -(27+7\Bi h) \varphi\right\}
		\\ && + J_\varphi \frac{(\P_x h)^2}{h^2} + \frac{4}{h} \P_x h \P_x \varphi
		+L_\varphi   \frac{\P_x h \P_x \theta}{h} + \P_{xx} \varphi
		\label{eq:model-theta-phi-varphi}
	\end{eqnarray}
	referred hereinafter as the $\theta$~-~$\varphi$~model,
	with
	\begin{eqnarray}
		\nonumber
		E_\varphi &=& - \frac{3(25+11\Bi h)}{14} \,,  \qquad  F_\varphi =- \frac{66+9\varphi+ 6(38\Bi h -11)\theta }{28}\,,\qquad G_\varphi= \frac{57}{7} \Bi h \,,\\
		J_\varphi &=& 6-(25+7\Bi h)\varphi + 6(2\Bi h -1)\theta\,,\qquad
		L_\varphi =48 -12 \Bi h -14 (\Bi h)^2\,.
		\label{eq:model-theta-phi-coeffs}
	\end{eqnarray}
\end{subequations}
The evolution equation \eqref{eq:model-ts} is the trace of the Fourier equation taken at the interface. As a consequence, it is exact and independent of the choice of the polynomials ${\hat v_1}$ and ${\hat v_2}$.

By construction, model~\eqref{eq:model-theta-phi} is consistent at order $\eps$. A study of the linear response of the model to a sinusoidal perturbation of wavenumber $k$ assuming a uniform film flow (i.e. $h=1$ and $q=1/3$) yields the matrix
\begin{equation}
	\label{eq:lpm}
	C = \begin{bmatrix}
		- k^2                   & 1                     \\
		-60 \, (1 + \Bi\,h) & -(27 + 7 \, \Bi \, h) - k^2
	\end{bmatrix}
\end{equation}
whose eigenvalues $\lambda_\pm$ are compared to the two first eigenvalues \eqref{eq:lin} and to damping rate \eqref{eq:l-theta} in figure~\ref{fig:lambda}. $\lambda_+$ is a good approximation to $\lambda_1$ whereas $\lambda_-$ is a poorer one to  $\lambda_2$. Nevertheless, $\lambda_+ \approx \lambda_1$ shall guarantee that the diffusive relaxation to the linear temperature distribution is correctly captured by the model.

\begin{figure}
	\begin{subfigure}[t]{0.49\textwidth}
		\includegraphics{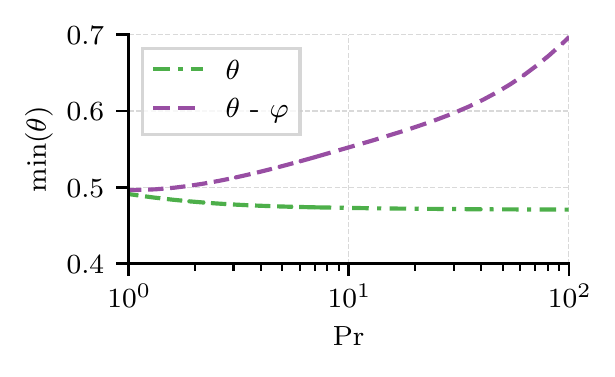}
		\caption{$\Bit=0.1$}
	\end{subfigure}
	\begin{subfigure}[t]{0.49\textwidth}
		\includegraphics{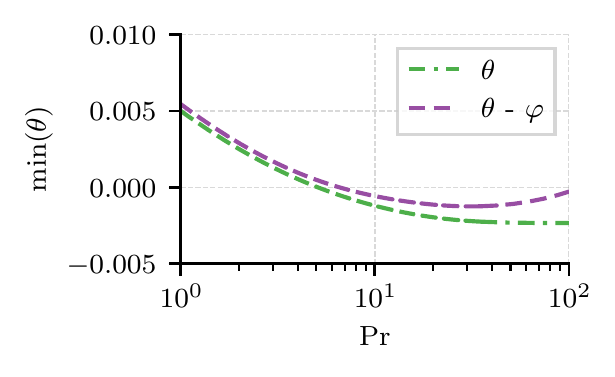}
		\caption{$\Bit=10$}
	\end{subfigure}
	\caption{Minimum value of a $\theta$ under a solitary wave as function of the Prandtl number ($\Re=15$, $\beta=90^\circ$ and $\Ka=3000$).
		\label{fig:mintheta}}
\end{figure}

Figure~\ref{fig:mintheta} compares the minimum values of the free surface temperature obtained with the one-variable \eqref{eq:model-theta} and two-variable model \eqref{eq:model-theta-phi} for two Biot number. A minor improvement is observed using two variables for a high Biot number instead of the one-variable model. For both models, $\min(\theta)$ presents nonphysical negative values. The two-variable model remains closer to $0$ than the $\theta$ model. For both models, this nonphysical behavior is limited compared to previous attempts.

\begin{figure}
	\begin{center}
		\begin{subfigure}[b]{0.49\textwidth}
			\includegraphics{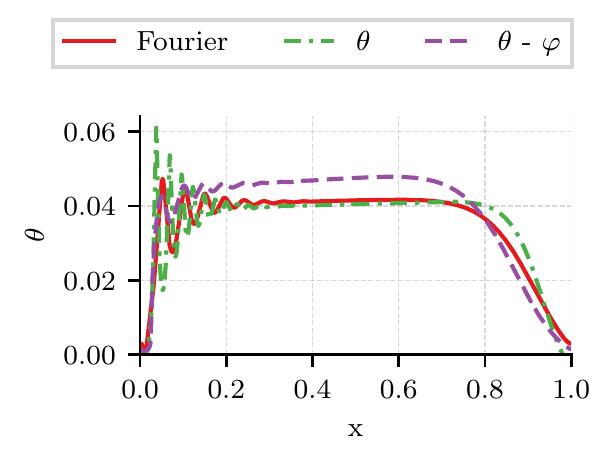}
			\caption{free surface temperature $\theta$ according to $x$}
			\label{fig:prof-NSF-phithCL-xtheta}
		\end{subfigure}
		\begin{subfigure}[b]{0.49\textwidth}
			\includegraphics{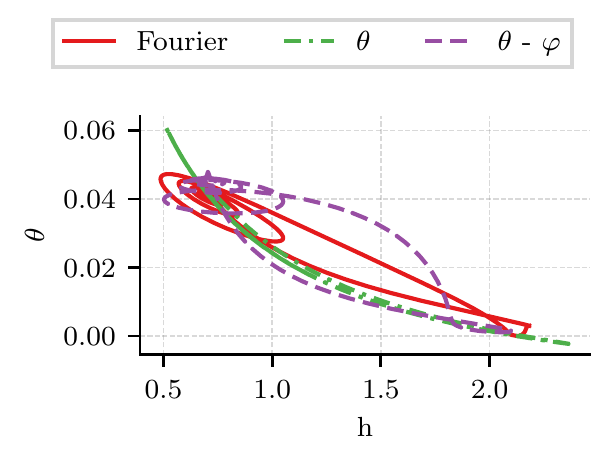}
			\caption{free surface temperature $\theta$ according to $h$}
			\label{fig:prof-NSF-phithCL-htheta}
		\end{subfigure}
	\end{center}
	\caption{Film of water on a vertical wall at $Re=15$, $\Bit=10$, $\Pr=100$ and $f=8$~Hz.}
	\label{fig:prof-NSF-phithCL}
\end{figure}

Figure~\ref{fig:theta_thetah} has been completed with the results of the $\theta$~-~$\varphi$~model \eqref{eq:model-theta-phi}. A very noticeable improvement over the $\theta$~model \eqref{eq:model-theta} can be observed as the sharp variation of the free surface temperature at the hyperbolic stagnation point in the wave is precisely captured by the $\theta$~-~$\varphi$~model. This agreement has been obtained over a wide range of $\Bi$ and $\Pr$ number.

Figure~\ref{fig:prof-NSF-phithCL} presents a comparison in the very demanding case of large values of $\Pr$ and $\Bi$ numbers where the agreement to the Fourier solution is the least convincing. Yet, the solution to the $\theta$ model agrees again well with the Fourier solution in the wave tail, where the film is close to the Nusselt solution (a flat film), but has some trouble to remain accurate as the shape of the wave becomes more complex. It is especially obvious when $\theta$ is plotted according to $h$ (fig. \ref{fig:prof-NSF-phithCL-htheta}). Adding a second variable, more of the temperature surface distribution complexity is captured, $\theta(x)$ mimicking well the reference solution. The simplest model is still advantageous : it shows a good accuracy to capture the averaged properties along the wave where the $\theta$~-~$\varphi$~model fail to predict the surface temperature where the film is almost flat. This is a common behavior of complex models: they improve accuracy and are able to capture more complex phenomena but are less robust and fail when the case is more demanding and far away from the asymptotic (here an order-one Peclet hypothesis). This can be observed for the long-wave Benney equations which capture the hydrodynamics of the film at low values of the Reynolds number.
The second-order Benney equation,  even if more accurate than the first-order one, is unable to deal with moderate Reynolds numbers \cite{gottlieb2004}. \label{review:complex}
\begin{figure}
	\begin{center}
		\includegraphics{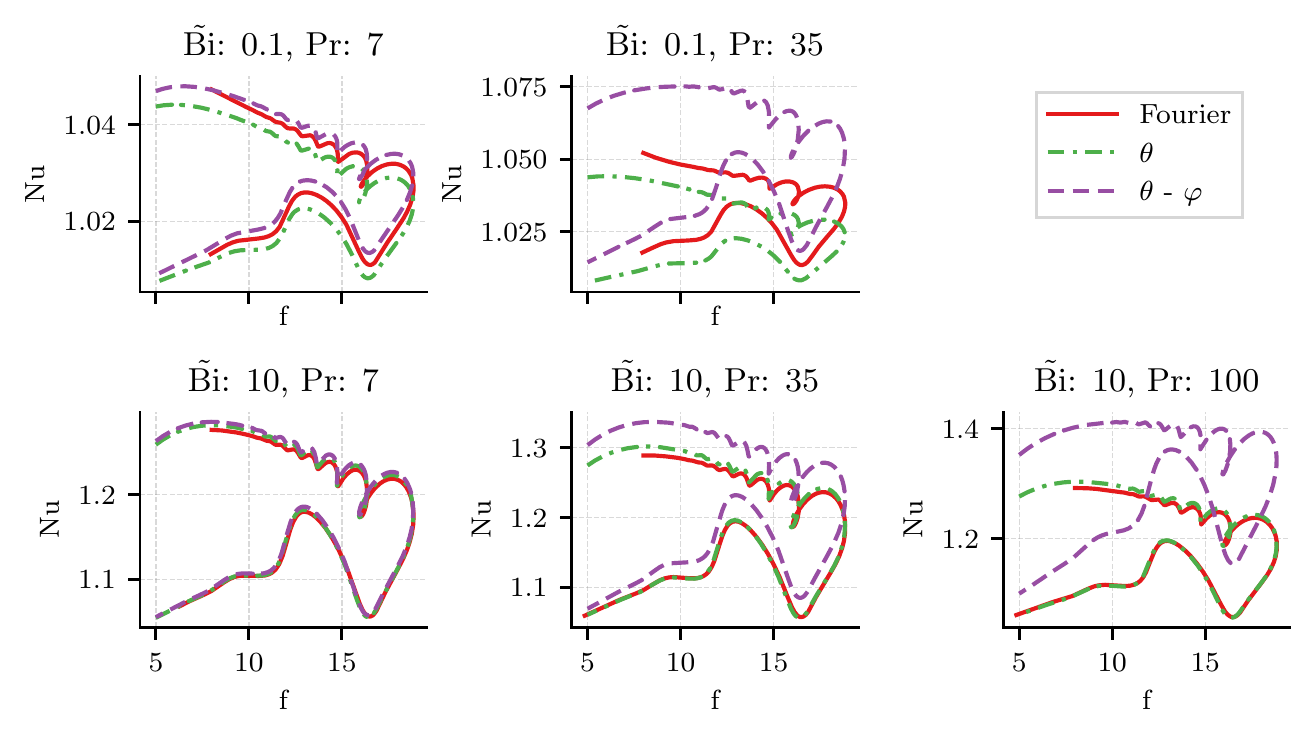}
	\end{center}
	\caption{Film of water on a vertical wall at $Re=15$. The graphs show the global Nusselt number (average of the flux at the free surface rescaled by its value for a flat film). Top $\Bit=0.1$, Bottom $\Bit=10$. First column $\Pr=7$, second column $\Pr=35$, last column $\Pr=100$.}
	\label{fig:Nu_f}
\end{figure}

This is confirmed by the results displayed in figure \ref{fig:Nu_f}, which presents the global Nusselt number for traveling-wave solutions (computed as the average of the flux at the free surface rescaled by its value for a flat film) according to the wave frequency. The flux being averaged, the $\theta$~-~$\varphi$ model main advantage (the ability to represent the complexity of the heat transfer in a more complex hydrodynamic regime) recedes, and the $\theta$ model performs somewhat better, especially for high Prandtl number. However both models capture accurately the global heat flux through wave in the thermally developed regime. This is particularly true dealing with water ($\Pr = 7$). Departures from the predictions of the Fourier equation can be observed at high values of Prandtl number. Yet, both models provide reasonable answers even at $\Pr = 100$.

To conclude, the two models \eqref{eq:model-theta} and \eqref{eq:model-theta-phi} have different advantages. The first one is robust, and can lead to a better global accuracy. It also uses only one variable to parametrize the thermal transfer, leading to cheaper resolution cost. The latter is able to represent more complex behaviors at a cost of a somewhat lower robustness (and global accuracy) and a higher computational cost (which is still far less expensive than solving the full Fourier equation). According to the goal of the study, one or the other may be used.

\section{Time dependent simulations}
\label{section:results}

The proposed formulations for the averaged heat balance have been validated through computations of the traveling-wave solutions, which implies a thermally and hydrodynamically developed regime. However, describing accurately the entrance region of a film flow developing on a plate is crucial for the optimization of a plate exchanger. Therefore, we turn to time-dependent simulations of heat transfer across a 2D liquid falling film. These simulations have been performed using the Saint-Venant hydrodynamic formulation, proposed by \citet{ruyer-quil2000}, which reads as

\begin{equation}
	\begin{aligned}
		\begin{cases}
			 & \partial_{t}h = -\partial_{x}q                              \\
			 & 3  \Re \partial_{t}q = \frac{5\,h}{6} - \frac{5\,q}{2\,h^2}
			\\ & \qquad
			+ 	\frac{3}{7} \Re \left(9 \partial_{x}h \frac{q}{h}
			-17 \partial_{x}q
			\right) \frac{q}{h}
			- \frac{5}{6} \Ct\,\partial_{x}h
			+ \frac{5}{6} \We\,h\,\partial_{xxx}h
			\\ & \qquad
			+\frac{4\,q}{h^2}\partial_{x}h^2
			- \frac{9}{2\,h}\partial_{x}h\partial_{x}q
			- \frac{6\,q}{h}\partial_{xx}h
			+ \frac{9}{2}\partial_{xx}q\,.
		\end{cases}
	\end{aligned}
	\label{eq:03:sv2}
\end{equation}

The reason of this choice is the model's capacity \eqref{eq:03:sv2} to adequately capture the nonlinear wavy regime of liquid falling films at low to moderate values of the Reynolds number, as demonstrated by comparisons to direct numerical simulation (DNS) (see for instance \citet{ruyer-quil2014}). In section \ref{section:primitive-equations}, the hydrodynamic model has been chosen to make a comparison with the previous study by \citet{chhay2017}.
The hydrodynamic parameters are \(\Re=15\), \(\Ct=0\) and \(\We=266\) in each case. They correspond to a water film flowing on a vertical plate. This relatively low value of the Reynolds number ensures that the hydrodynamics of the film is adequately captured by the model.

In parallel to the models of heat transfers, we also solved the basic Fourier equation \eqref{eq:fourier} to provide means of validations. To solve the Fourier equation, a change of coordinates has been performed with \(\bar{y}=y/h \in [0, 1]\) instead of \(y \in [0, h]\). As a result, the numerical domain is a fixed rectangle \(x\in[0, L]\), \(\bar{y}\in[0, 1]\), removing the need of a moving mesh.

\subsection{Application case example}

As an introduction to comparisons between models and their validation to the Fourier basic equation, a typical case is presented here. It corresponds to a water falling film flowing on a vertical plate with an oscillation of its inlet fluid height, periodic in time. We chose a moderate Reynolds number and a low Biot number (as seen in heat exchangers). The parameter set is the following:  $\Ct=0$, $\Re = 15$, $\We = 266$, $\Pe=105$, $\Bi = 0.1$.
The plate is maintained at a hot constant temperature, whereas the fluid flows in contact with a cold atmosphere.

\begin{figure}[btp]
	\centering
	\begin{subfigure}[b]{\textwidth}
		\includegraphics{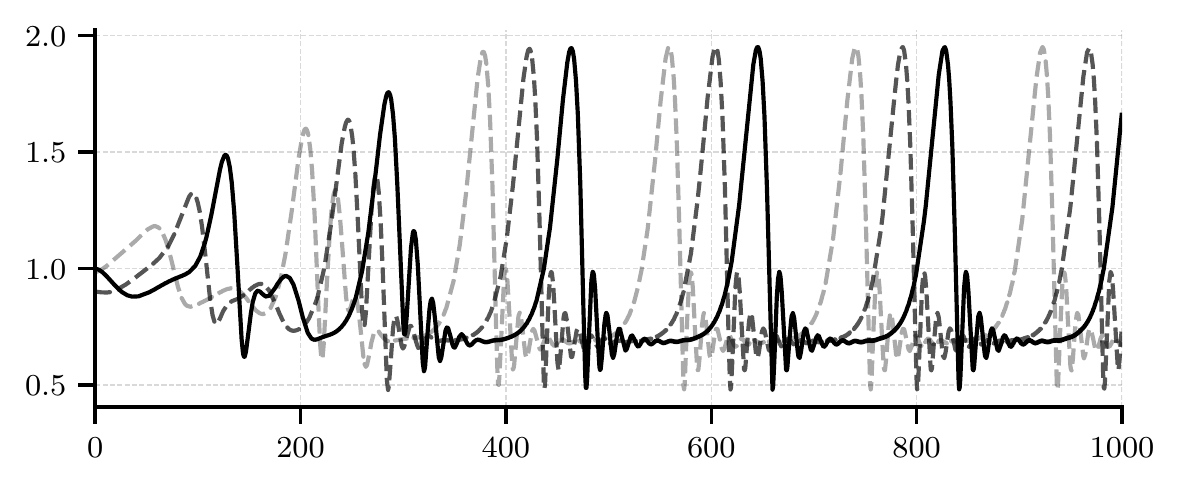}
		\caption{Film height for three successive snapshots.}
		\label{fig:film_height}
	\end{subfigure}

	\begin{subfigure}[b]{\textwidth}
		\includegraphics{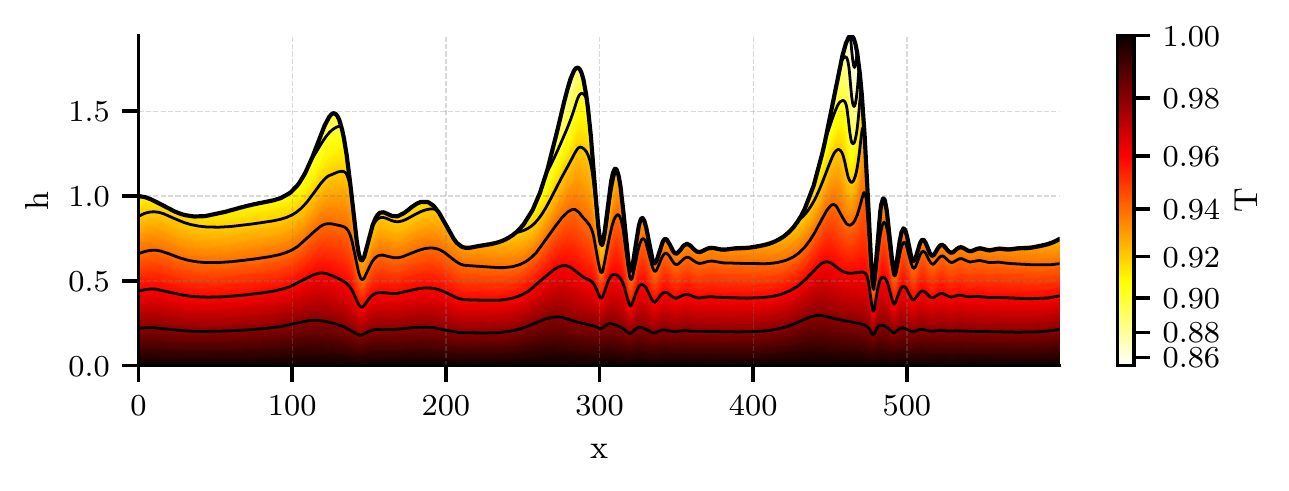}
		\caption{Snapshot of the fluid temperature.}
		\label{fig:film_h_T}
	\end{subfigure}
	\caption{Example study, water flowing over a vertical plate. $\Ct=0$, $\Re = 15$, $\We = 266$, $\Pe=105$, $\Bi = 0.1$. Temperature field computed solving the Fourier equation.}
	\label{fig:example_case}
\end{figure}

As we can see in figure \ref{fig:film_height}, the inlet oscillations grow quickly, leading to a saturated wavetrain. These waves consist of one main hump preceded by capillary waves. These capillary waves are close to each other and have a smaller amplitude than the main hump. Without inlet noise, these waves are evenly spaced and stable in time. Figure \ref{fig:film_h_T} shows the temperature field across the film. With moderate-to-high Peclet numbers, we can notice convective effects at the top of the main crest where the cold fluid near the interface mixes slightly with the fluid in the bulk region of the film.

\subsection{Comparison between models}
\label{subsec:comparison-between-the-models}

Simulations have been first run for a low Peclet number, in order to check the coherence with the Fourier equation. In the low Peclet limit assumption, where the long-wave expansion holds,
the temperature fields predicted by the model should agree with the reference solution to the Fourier equation.
\begin{figure}
	\centering
	\includegraphics{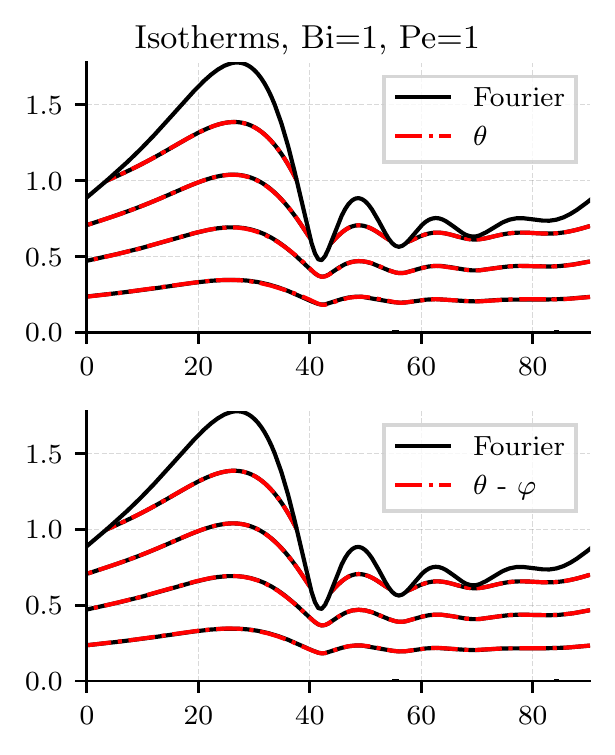}\hfill
	\includegraphics{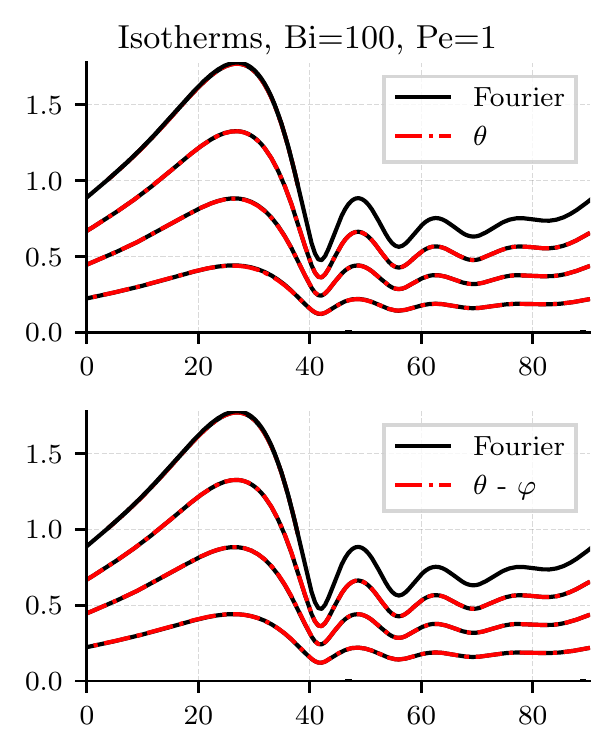}
	\caption{Comparison between the reference case and the different models at low Peclet ($\Pe = 1$). Left, $\Bi=1$ and right, $\Bi=100$. From top to bottom, we have the $\theta$~model and the $\theta$~-~$\varphi$ with polynomial test functions.}
	\label{fig:check_asympt}
\end{figure}
As we can see in figure \ref{fig:check_asympt}, in the limit case $\Pe \rightarrow 0$, the models present the same behavior as the reference Fourier model, for both moderate and high Biot numbers.

\begin{figure}
	\centering
	\includegraphics{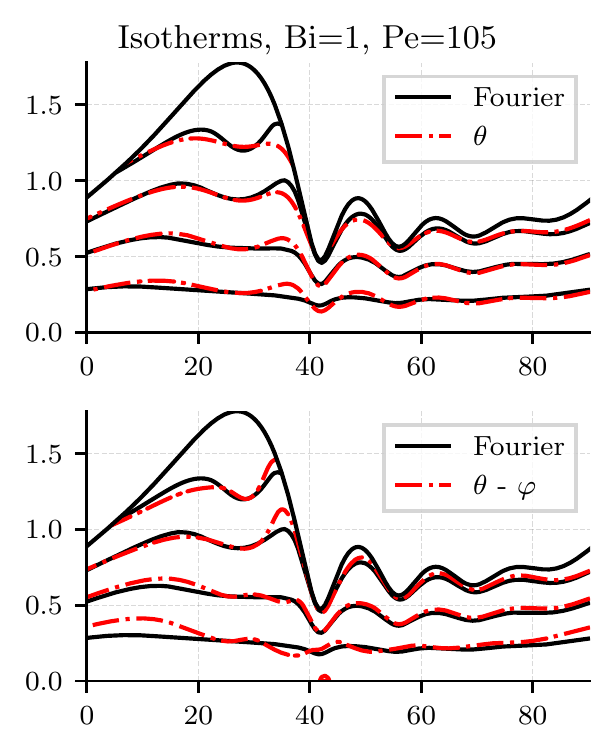}\hfill
	\includegraphics{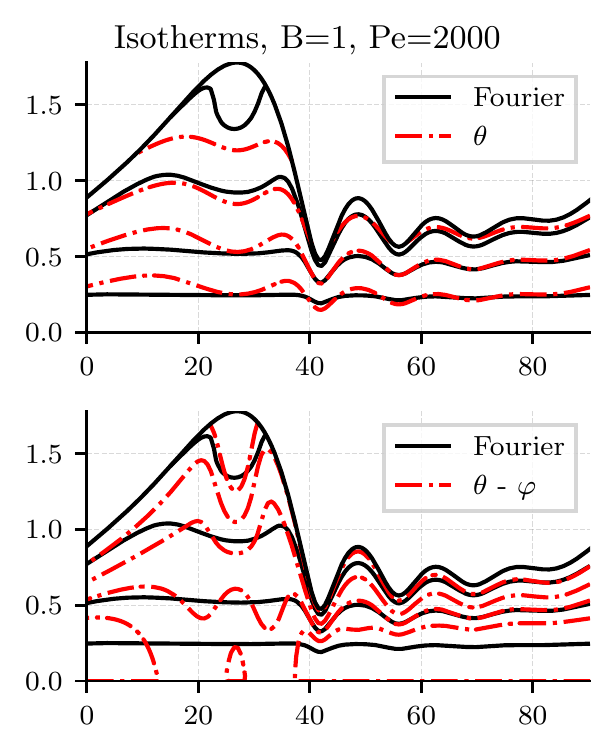}
	\caption{Comparisons between the reference case and the different models at moderate Biot ($\Bi = 1$). Left panel: high Peclet number $\Pe=105$ (coherent with water thermal properties). Right panel: very high Peclet number $\Pe=2000$.}
	\label{fig:high_peclet}
\end{figure}

As we increase the Peclet number, we still observe a good agreement with the Fourier equation, even if our models are built on a low Peclet hypothesis (cf. figure \ref{fig:high_peclet}). As stated previously, the $\theta$~model is not complex enough to catch the detail of the temperature field (especially in the crest, near the thermal boundary onset) where the $\theta$~-~$\varphi$~models are capable to exhibit a similar complexity. Considering the isotherms close to the wall, the simplest model shows a better agreement with the Fourier equation than the more complex one, where some spurious oscillations can be seen. This is not surprising, as the model is parameterized with only one variable corresponding to the liquid-gas interfacial temperature (where the temperature field presents the greatest complexity).
In any case, considering the fact that our derivation assumes order-one values of the Peclet number, the two models show a good agreement with the Fourier equation.

In addition to the periodic-box simulation, simulations of the evolution of the film in an open large domain, representative of an element of a plate exchanger, have been performed. This is a critical case, as our previous attempts \cite{cellier2018b} were unable to capture correctly the onset of the thermally developed regime at the inlet of the flow. The reason for this inaccuracy lies in an incorrect representation of the diffusion relaxation modes discussed in the previous section.
The thermal entrance length increases with the Peclet number and can exceed the exchanger length: this is an important factor for the heat exchanger optimization.
The same parameters as the periodic-box case have been chosen (\(\Re=15\), \(\Ct=0\), \(\We=266\)). We modeled a \(L=\unit{20}{\centi\meter}\) length exchanger plate. A Dirichlet boundary condition has been used at the flow input such as

\[h|_{x=0} = 1 + A\,\sin(2 \pi\,t\,f) \quad q|_{x=0} = \frac{1}{3}h^3\]

with the amplitude \(A=0.1\) and the frequency \(f=10\).

The outlet is dealt with a no-flux boundary condition: it yields some numerical errors that are convected outside the domain. We lose a small part of the simulation domain length, and therefore extend the domain to \(L=\unit{25}{\centi\meter}\). We then crop a buffer zone to obtain 20 cm of useful length for the simulation.

\begin{figure}[h!tbp]
	\centering
	\includegraphics{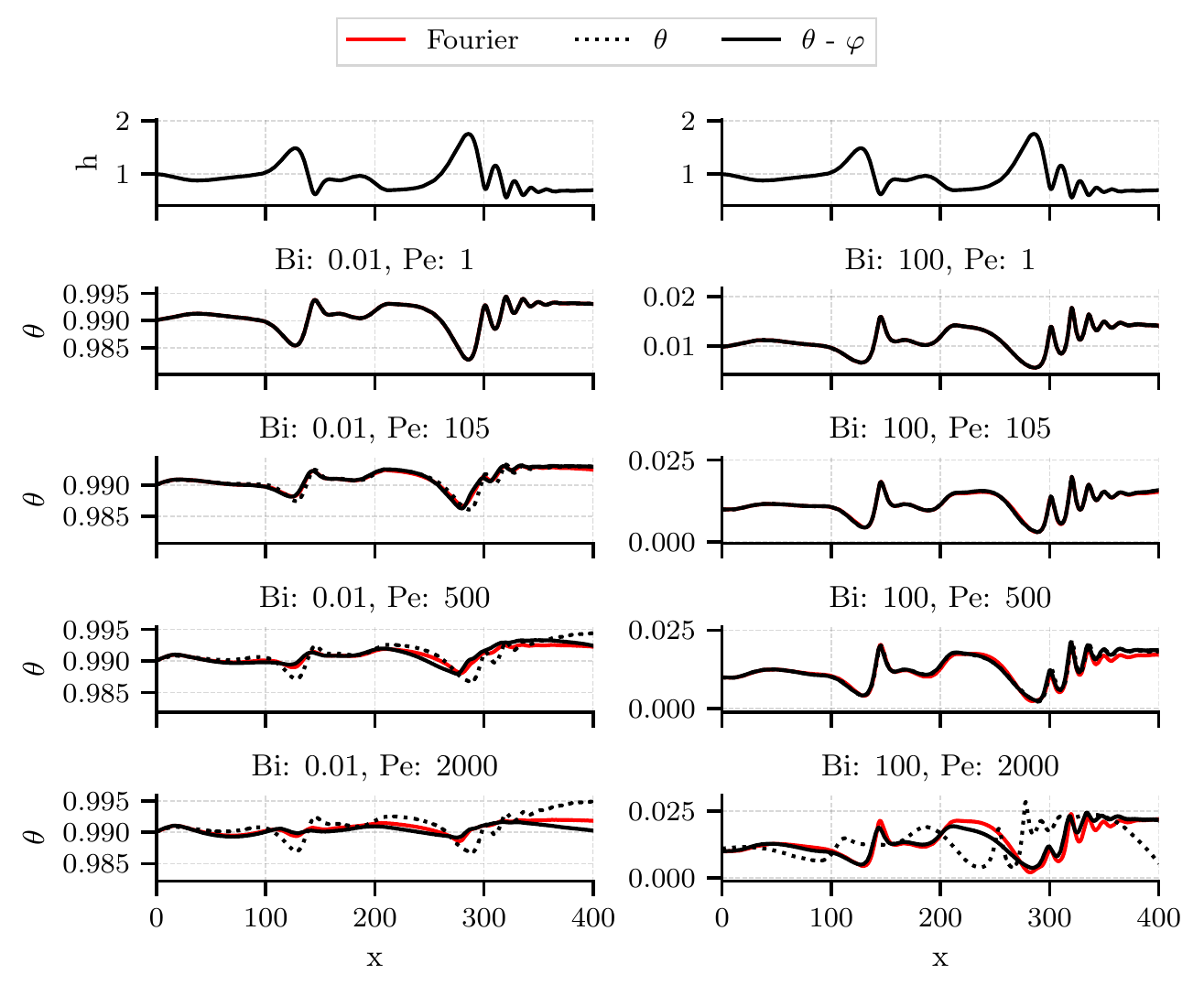}
	\caption{Simulation of a full exchanger plate, first part of the domain. \(\Re=15\), \(\Ct=0\), \(\We=266\).}

	\label{fig:start_plate}
\end{figure}

Figure \ref{fig:start_plate} focuses on the first part of the plate, where the waves are growing. The $\theta$~-~$\varphi$~model has the same behavior when the Biot number is low and shows a slightly better agreement with the Fourier equation than the $\theta$~model. Both over-estimate the interfacial temperature.

\begin{figure}[hb]
	\centering
	\includegraphics[width=\textwidth]{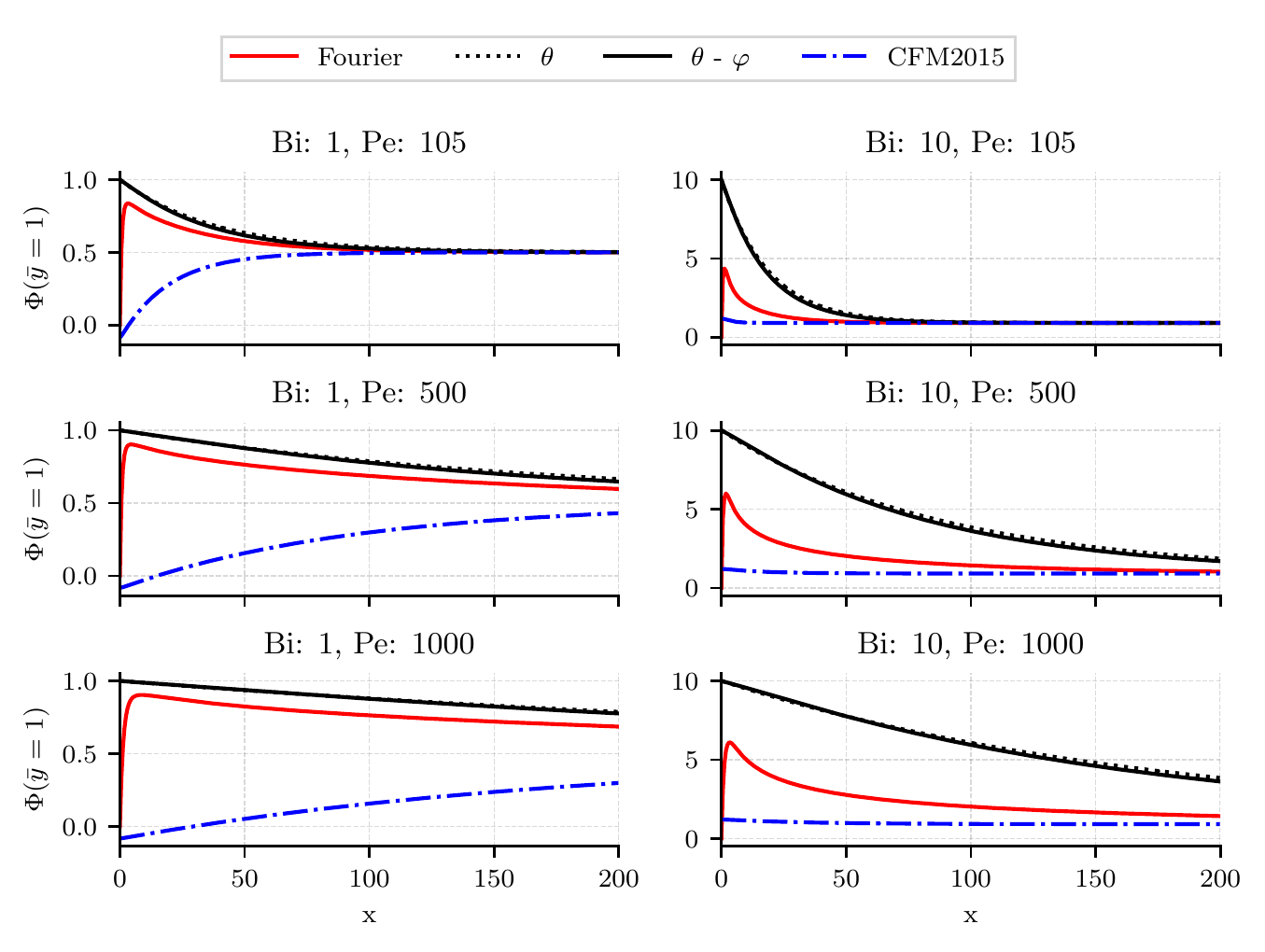}
	\caption{Interfacial heat flux along the flow length, flat film. \(\Re=15\), \(\Ct=0\), \(\We=266\).}
	\label{fig:relax_plate}
\end{figure}

To check the accuracy of the models for a relaxation process, some simulations have been run for a flat case (without any film perturbation) and a hot film input ($T|_{x=0}$ = 1). For the interfacial flux, our two new models (\(\theta\) and \(\theta-\varphi\)) have very close behaviors (see figure \ref{fig:relax_plate}). We are not able to capture the very first part of the relaxation, where the Fourier model goes from no flux at all to local maxima before relaxing. Our model cannot capture such a sharp transition, as a polynomial projection of the temperature field cannot represent a Dirac function. That explains the observed initial flux overshoot. A previous attempt (referred as CFM2015 \cite{ruyer-quil2015b})\label{review:CFM} is unable to capture the relaxation of the interfacial flux at all.

\begin{figure}[hb]
	\centering
	\includegraphics[width=\textwidth]{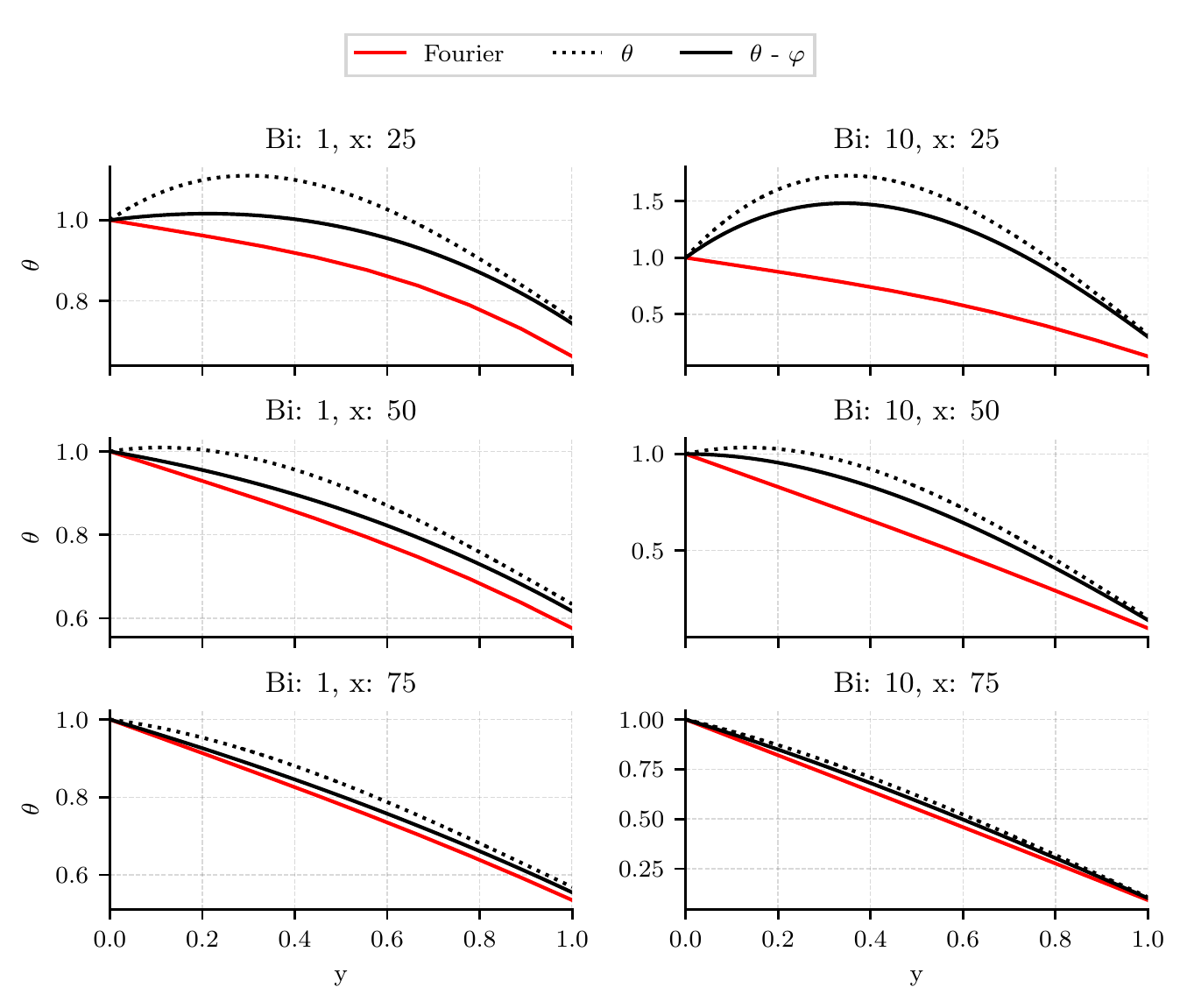}
	\caption{Temperature field profile at different position for a flat film. \(\Re=15\), \(\Ct=0\), \(\We=266\), \(\Pe=105\).}
	\label{fig:flat_temperature_profile}
\end{figure}

Similarly to \citet{aktershev2017}, we have represented the temperature profile of the different models as well as the solution of the Fourier equation for a flat film. This has been done at different positions on the plate, as shown in figure~\ref{fig:flat_temperature_profile}. As stated previously, the polynomial projection of the temperature field, parametrized with interfacial-based free variable, cannot capture a sharp transition. This led to an unphysical representation near the plate (especially for the $\theta$~model). Yet, the interfacial temperature is correctly captured. The position where the linear temperature profile is reached is predicted more accurately when the Biot number is low. This is a consequence of our polynomial approximation which is linear with respect to the low Biot number and can be improved by using more test functions and refining the projection. This improvement will lead to higher model complexity.

\subsection{Validation - periodic box}

A series of simulations have been computed with fixed hydrodynamic parameters, the only varying parameters being the Biot and the Peclet number. The different simulations are compared to the Fourier reference case.

For very large Peclet numbers, this procedure is not sufficient to guarantee an accurate representation of the temperature field, especially in the vicinity of the thermal boundary layer. However, the obtained accuracy is adequate for the validation of the models.

\begin{figure}
	\centering
	\begin{subfigure}[t]{0.47\textwidth}
		\includegraphics{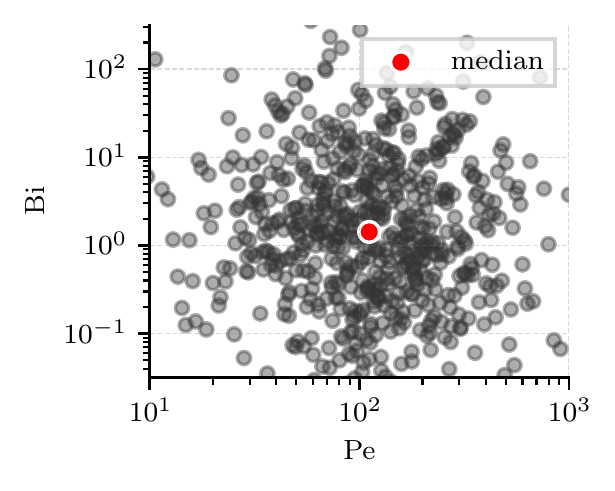}
		\caption{Periodic box study case : 640 samples.}\label{fig:validation_sample}
	\end{subfigure}
	\hfill
	\begin{subfigure}[t]{0.47\textwidth}
		\includegraphics{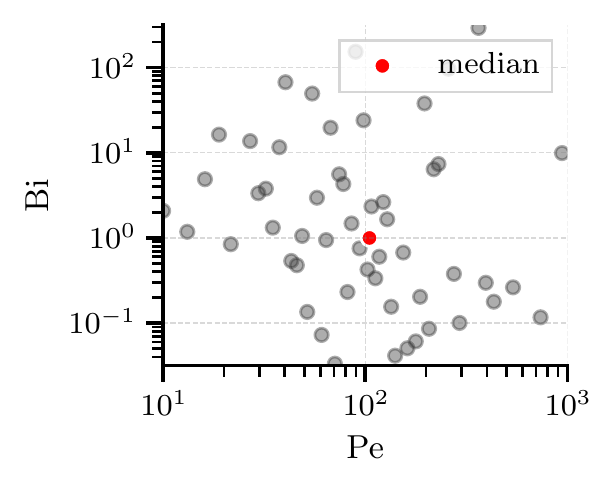}
		\caption{Open flow study case : 64 samples.}\label{fig:open_validation_sample}
	\end{subfigure}
	\caption{Validation sampling: samples chosen with the latin hypercube sampling (LHS)  method.}
\end{figure}

The chosen sampler is a Latin Hypercube Sampler~\cite{mckay1979} generating samples following a log-normal distribution. The log-normal shapes are chosen in order to fix the median for both varying parameters. The samples are summed up in the figure~\ref{fig:validation_sample}.
The number of samples (640) is large enough to provide a good overview of the behavior of the models according to the two varying parameters.

The two models presented in the previous section are used to simulate a traveling wave in a periodic box of length \(L/\bar{h_N} = 90\). The long-time solution of the different models is compared with the reference solution to the Fourier equation.
\begin{figure*}
	\centering
	\includegraphics[width=\textwidth]{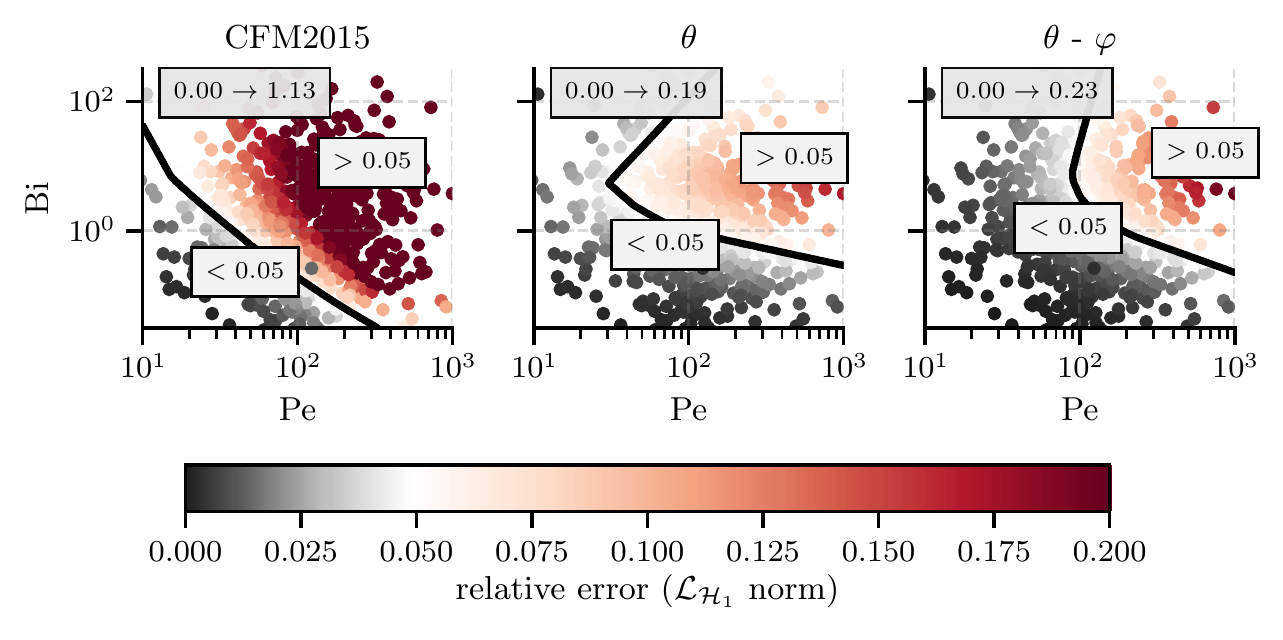}
	\caption{Relative error (\(\mathcal{H}_1\) norm) on the wall and interface flux, for the different models. Minimum and maximum error values are displayed at the top left of each plots, and the black border separates the domain where the error is inferior to 5\%. This border is determined by training a multi-layered perceptron (MLP) classifier with our data.}\label{fig:errorH2}
\end{figure*}
Figure~\ref{fig:errorH2} compares the different results. The norm \(\mathcal{H}_1\) is defined as
\begin{equation}
	\mathcal{L}_{\mathcal{H}_1}(X) = \sqrt{\int{X^2 + \frac{\partial X}{\partial x}^2 \mathrm{d}x}}\,.
	\label{eq:LH2}
\end{equation}
This norm has been chosen to evaluate both the amplitude and the shape of the heat flux distributions at the interface as predicted by the models.
We can observe that the two models exhibit a very good agreement with the reference solutions. For more than half of the physical domain investigated, the error is below a 5\% threshold, and never exceed 25\%. The $\theta$~-~$\varphi$~model presents a wider domain in the parameter space of applicability, defined here by the 5\% error threshold. In particular, the range of validity for the $\Pe$ number has been extended up to $\Pe = 100$, whatever the value of the $\Bi$ number, which is a significant improvement over the $\theta$~model. However, the $\theta$~model presents less pronounced maxima of deviation from the Fourier solutions (with a maximal error of $19\%$ instead of $23\%$) but is not able to represent some important phenomena, such as the developing thermal layers near the crest of the waves (as shown in figure \ref{fig:high_peclet}).

\subsection{Validation - full exchanger}

The same set of parameters as the periodic-box case has been chosen for the validation case. We have simulated a \(L=\unit{20}{\centi\meter}\) length exchanger plate. The same boundary conditions as in section \ref{subsec:comparison-between-the-models} have been used. A regular forcing at the inlet is again enforced with an amplitude equal to \(A=0.1\) and a frequency given by \(f=10\).

Such simulations being expensive in comparison with the periodic-box case, we limit the sample number to 64. The samples can be seen in figure \ref{fig:open_validation_sample}: the parameter space is well explored and the shape of the log-normal distribution has been chosen so that the median is aligned with our reference case.

\begin{figure*}
	\centering
	\includegraphics[width=\textwidth]{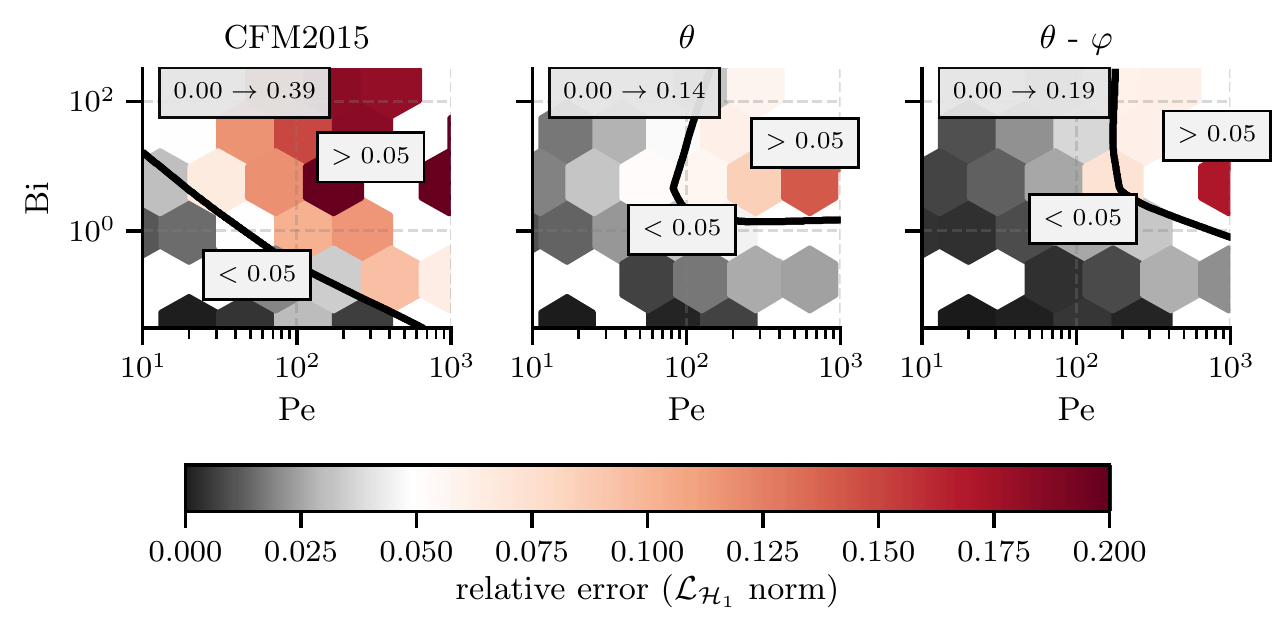}
	\caption{Relative error (\(\mathcal{H}_1\) norm) on the wall and interface flux, for the different models. The minimum and maximum error values are displayed at the top left of each plot, and the black border separates the domain where the error is inferior to 5\%. This border is determined by training a MLP classifier with our data.}
	\label{fig:open_errorH2}
\end{figure*}

As observed in figure \ref{fig:open_errorH2} (see \eqref{eq:LH2} for the \(\mathcal{H}_1\) norm definition), the error of the transient state is smaller than the error for a steady traveling wave, and we have seen that our models relax well to the equilibrium state. The same remarks made for the periodic box stay: the $\theta$~model fails to represent the complexity of the temperature field (see figure \ref{fig:high_peclet}). Yet, this simple model catches well the interfacial flux with an error below 15\%, which is a strong improvement compared to the previous attempts \cite{cellier2017,chhay2017}.

\section*{Conclusion}

A new asymptotic model, offered with two variants, has been developed as an alternative to the full resolution of the Fourier equation across a falling film. This leads to a faster resolution at the cost of a reduced domain in the parameter space of applicability (very high \Pe numbers are still out of reach). This model overcomes the limitations of previous attempts, which led to acceptable results for moderate Biot and Peclet numbers only, and yielded non-physical behavior outside this range. Moreover, the diffusive relaxation towards the conductive equilibrium in the entrance region of the plate observed with the Fourier equation is now correctly captured. This improves the models accuracy, even outside the entrance regime. It also extends the physical space of applicability (see figure \ref{fig:errorH2}, \ref{fig:open_errorH2}), even with only one free variable to represent the temperature distribution (relative error less than 20\% for $\Pe \in [10^1, 10^3]$ and $\Bi \in [10^{-3}, 10^3]$).

The two variants, resp. $\theta$~model and $\theta$~-~$\varphi$~model, possess different advantages. The simplest one ($\theta$~model) is more robust and has a cheaper resolution cost. It is a good candidate for global studies (optimization for example), whereas the second one ($\theta$~-~$\varphi$~model) is able to capture more complex thermal transfer behaviors. This complexity has a cost, in terms of robustness and computation especially, at the crest of the waves (due to the extra free variable and evolution equation). That computation cost is still much less expensive than the alternative, i.e. solving the Fourier equation. Hence, this model shows itself to be a good candidate when the comprehension of the phenomena is important but the cost of the Fourier equation cannot be afforded. This cost can be prohibitive when it comes to transfers within a 3D falling film in a spatial domain representative of a realistic plate exchanger. The latter is our next goal, within reach by coupling our model with a computation-efficient shallow-water model developed recently by \citet{bresch2019}. Other perspectives include the introduction of coupling effects between hydrodynamic and heat transfer via the Marangoni effect, or via other temperature dependencies of the fluid properties \cite{pascal2019} in the models. The introduction of such coupling is trivial and will give access (with the extension to 3D of the models) to a proper comparison with the experimental studies. \label{review-perspectives} \tagreview{Readers interested in how to extend such models in 3D or how to account for the thermocapillarity can find details in \citet{kalliadasis2012}. \label{review-book}}

This family of models constitutes a new tool which provides a costless evaluation of thermal transfers across a falling film, making costly investigations in terms of number of simulations (optimization, sensitivity analysis, parameterized exploration of the parameter space...) now accessible.

\appendix

\section{Construction of solutions to the Fourier equation}
\label{sec:appendix2}
In this section, we present an approach to solve the Fourier equation \eqref{eq:fourier} in the case of traveling wave solutions. 
We first project the temperature field as
\begin{equation}
\label{eq:Cheb}
    T = 1 +  \sum_{i=1}^n \tau_i(x,t)\phi_i(X)\quad \hbox{and}\quad
    X = 2\frac{y}{h(x,t)} -1
\end{equation}
where $\phi_i(X)$ are linear combinations of Chebyshev polynomials of the first kind $T_i$ given by
\begin{eqnarray}
\nonumber
\phi_1 &=& 1 + X \,,\\
 \phi_{2i} &=& T_{2i} (X) -1 \quad \hbox{and}\quad
 \phi_{2i+1} = T_{2i+1} (X) -X
 \quad \hbox{for} \quad i \ge 1\,,
\end{eqnarray}
so that $\phi_i(1)= \phi_i(-1)=0$ for $i\ge 2$ and $\phi_1(-1)=0$. The Dirichlet condition at wall \eqref{eq:CMA-BC-1} is thus verified by \eqref{eq:Cheb}. Considering traveling waves, i.e. stationary solutions in frame
$\xi=x - c\, t$ moving at a constant speed $c$, and writing the
Fourier equation \eqref{eq:fourier}  on the Gauss-Lobato points 
$X_i = - \cos(\pi i/n)$ $i\ge 1$ gives formally $n-1$ relations 
\begin{equation}
\label{eq:Fourier-Cheb}
    \sum_j^n \phi_j (X_i) D_{\xi\xi} \tau_j = F_i(\tau_j, D_\xi \tau_j)\,,
\end{equation}
where $D_\xi = d/d\xi$.
We next complete the Fourier equation \eqref{eq:fourier} by the boundary condition
\begin{equation}
\label{eq:Newton-eta}
    \eta \partial_{xx} T = \P_y T - \P_x h \P_x T
		  +  \Bi T \sqrt{1 + (\P_x h)^2} \qquad\hbox{at}\quad y=h\,.
\end{equation}
 The Newton law of cooling \eqref{eq:CMA-BC-2} is recovered in the limit $\eta\to0$.
Substitution of \eqref{eq:Cheb} into \eqref{eq:Newton-eta} completes the $n-1$ relations \eqref{eq:Fourier-Cheb} into a linear system of dimension $2n$
\begin{equation}
\label{eq:Fourier-system}
{\bf A} \frac{d {\bf U}}{d\xi} = {\bf B}({\bf U}; \eta)\,,
\end{equation}
with ${\bf U} = (\tau_i, D_\xi \tau_i)$, $1\le i \le n$. Inverting \eqref{eq:Fourier-system} leads to
an autonomous dynamical system of dimension $2n$. This dynamical system is solved along with the dimension-three dynamical system corresponding to \eqref{eq:Vila} or \eqref{eq:03:sv2} with the help of the software {\sc AUTO07p} \cite{doedel2007} (see \cite{kalliadasis2012} for detail). The value of the constant $\eta$ has been set to $10^{-6}$. We checked the convergence with respect to $\eta$ by dividing its value by $10$.

\section*{Acknowledgements}
The authors acknowledge support by the FRAISE project, grant ANR-16-CE06-0011 of the French National Research Agency (ANR) and by the project Optiwind through Horizon 2020/Clean Sky2 (call H2020-CS2-CFP06-2017-01) with Saint-Gobain.
\FloatBarrier

\bibliographystyle{abbrvnat}
\bibliography{thermal_falling_film}

\end{document}